\documentclass[12pt]{amsart}
\usepackage{amsmath,amssymb,amscd,array, mathrsfs }

\usepackage{amsmath,amscd,amssymb,amsthm,array}
\usepackage{color}

\usepackage{amsmath,amssymb,mathrsfs,amsthm, tikz-cd,mathrsfs}
\usepackage{comment}
\def\url#1{\expandafter\s

\tring\csname #1\endcsname}

\def\mmat #1,#2,#3,#4,{\text{\small\arraycolsep=3pt $
\begin{pmatrix}#1&#2\\#3&#4\end{pmatrix}$}}

\usepackage{hyperref}

\usepackage{comments}
\newComments\SBe{Said}{blue}
\newComments\SBo{Sofiane}{blue}
\newComments\AM{Nacer}{blue}
\newComments\DL{DL}{red}
\newComments\QEh{QEh}{blue}

\def\mmat #1,#2,#3,#4,{\text{\small\arraycolsep=3pt $
\begin{pmatrix}#1&#2\\#3&#4\end{pmatrix}$}}

\usepackage{lscape}
\usepackage{tikz-cd}

\usepackage{DLdef1}
\usepackage{multicol}

\def\mmat #1,#2,#3,#4,{\text{\small\arraycolsep=3pt $
\begin{pmatrix}#1&#2\\#3&#4\end{pmatrix}$}}

\renewcommand {\ssbegin}[2][*]
 {\refstepcounter{subsection}%
\if#1*
\addcontentsline{toc}{subsection}{\thesubsection.\hskip 1pc #2}%
\else
\addcontentsline{toc}{subsection}{\thesubsection.\hskip 1pc #2. #1}%
\fi
 \def \secno {\gdef \secno {}{\ssecfont
\thesubsection.\hskip 2ex}%
 }%
 \begin{#2}}

\renewcommand {\sssbegin}[2][*]
 {\refstepcounter{subsubsection}
\if#1*
\addcontentsline{toc}{subsubsection}{\thesubsubsection.\hskip 1pc #2}%
\else
\addcontentsline{toc}{subsubsection}{\thesubsubsection.\hskip 1pc #2. #1}
\fi
 \def \secno {\gdef \secno {}{\ssecfont \thesubsubsection.\hskip 2ex}%
 }%
 \begin{#2}}

\renewcommand {\parbegin}[2][*]
 {\refstepcounter{paragraph}
\if#1*
\addcontentsline{toc}{paragraph}{\theparagraph.\hskip 1pc #2}%
\else
\addcontentsline{toc}{paragraph}{\theparagraph.\hskip 1pc #2. #1}
\fi
 \def \secno {\gdef \secno {}{\ssecfont \theparagraph.\hskip 2ex}%
 }%
 \begin{#2}}

\newcommand {\ce}{{\text{CE}}}

\rmnameii{etr}{etr}
\rmnameii{evv}{ev}
\DeclareMathOperator{\h}{\mathcal{H}}
\DeclareMathOperator{\K}{\mathbb{K}}

\newcommand{\black}{\color{black}}
\newcommand{\Z}{\mathbb{Z}}

\begin{document}

\title[Central extensions of restricted Lie superalgebras and classification]{Central extensions of restricted Lie superalgebras and classification of $p$-nilpotent Lie superalgebras in dimension $4$}

\author{Sofiane Bouarroudj}

\address{Division of Science and Mathematics, New York University Abu Dhabi, P.O. Box 129188, Abu Dhabi, United Arab Emirates.}
\email{sofiane.bouarroudj@nyu.edu}

\author{Quentin Ehret}
\address {Division of Science and Mathematics, New York University Abu Dhabi, P.O. Box 129188, Abu Dhabi, United Arab Emirates.}
\email{qe209@nyu.edu}

\thanks{The authors were supported by the grant NYUAD-065.}






\begin{abstract}
The first main result of this paper is to build the first and second restricted cohomology groups for restricted Lie superalgebras in characteristic $p\geq3$, modifying a construction by Yuan, Chen and Cao. We will explain how these groups capture some algebraic structures, such as extensions and derivations.  Further, we apply this construction to classify $p$-nilpotent restricted Lie superalgebras up to dimension $4$ over an algebraically closed field  of characteristic $p\geq3$.
\end{abstract}


\maketitle

\begin{flushright}
\textit{In memory of Pavel Grozman.}
\end{flushright}

\thispagestyle{empty}
\setcounter{tocdepth}{2}
\tableofcontents

\section{Introduction} \label{SecDef}

\subsection{Restricted Lie (super)algebras and their cohomology} Over a field of positive characteristic $p>0$, restricted Lie algebras are of prime interest, mainly due to their link to algebraic groups and their role in representation theory and classification (see \cite{SF}). Their definition goes back to Jacobson (\cite{J}), who derived it from the study of the space of derivations of a Lie algebra. A given Lie algebra $L$ is called \textit{restricted} if it is endowed with a so-called $p$-map $(\cdot)^{[p]}:L\rightarrow L$ that behaves roughly like the $p^{\text{th}}$ power in associative algebras, see Definition \ref{defres}. The cohomology associated with restricted Lie algebras is considerably more complicated than the ordinary Chevalley-Eilenberg cohomology (see \cite{CE}). In \cite{H,H1,H2}, Hochschild defined the restricted cohomology of a restricted Lie algebra $L$ with values in a module $M$ by
$$ H^k_{\rm{res}} (L;M) := \text{Ext}^k_{U_p(L)}(\K;M),~k\geq0,$$
where $U_p(L)$ denotes the restricted enveloping algebra of $L$. Although correct, this expression only allows explicit calculations for $k\in\{0,1\}$, in the context of certain extensions (see \cite{H}). Evans and Fuchs then proposed an explicit construction (for $p\geq 3$) of a cochain complex in \cite{E, EF}, which allows restricted cohomology groups to be computed up to order $2$ in the general case. Our understanding of this restricted cohomology is still incomplete, although the work of Evans and Fuchs has led to good cohomological interpretations of certain algebraic phenomena, see \cite{EF}. This cohomology was used to study central extensions of several classes of restricted Lie algebras; for instance, restricted Witt algebras (see \cite{EFP}), restricted filiform algebras (see \cite{EF2}) and restricted affine nilpotent algebras (see \cite{EF3}). Moreover, formal deformations were investigated in \cite{EM}. The superization of the notion of restricted Lie algebra is due to Mikhalev (\cite{Mi}), who called them $p$-superalgebras. Several authors have since studied this notion in various contexts. For instance, identities on the restricted enveloping superalgebras were investigated in \cite{P,U}, double extensions in \cite{BBH, BEM}, descriptions of $p|2p$-maps on simple restricted Lie superalgebras for $p=3$ in \cite{BKLLS}, and representations of the restricted Witt superalgebras and $\text{sl}(n|1)$ in \cite{SZ,Z}, respectively.
A generalization of restricted cohomology to restricted Lie superalgebras was investigated in \cite{YCC}. However, their construction only involves even elements, so it is not suitable to capture central extensions. Therefore, some results, like \cite[Theorem 3.7]{YCC}, have to be amended, see Theorem \ref{thmext}. The purpose of this paper is to propose a new definition for restricted cohomology of restricted Lie superalgebras (see Section \ref{restcohomosuper}) and investigate some classical interpretations (see mainly Section \ref{sectionextensions}). This paper does not consider any versions of restricted Lie superalgebras for $p=2$ as well as the corresponding restricted cohomology.

\subsection{Classification of $p$-nilpotent restricted Lie (super)algebras} A classical question related to the structure theory of Lie algebras is the classification problem. The application of root systems to the classification of semi-simple Lie algebras has been known since Cartan and Killing. In case that the Lie algebra is nilpotent, then the possibility of obtaining an explicit classification (by explicitly giving the multiplication table based on the basis elements) is possible. The classification of nilpotent Lie algebras over any field has existed for a long time up to dimension 5. Over a field of characteristic $0$ and dimension $6$, the classification was achieved by Morozov in 1958 (\cite{M}). In 2012,  a complete classification of nilpotent Lie algebras of dimension $6$ over an arbitrary field was achieved in (\cite{CGS,DG}). In the case where the Lie algebra is restricted, Schneider and Usefi have investigated $p$-nilpotent restricted  Lie algebra of dimension up to $4$ over a perfect field (including the case $p=2$) in \cite{SU}, then Darijani and Usefi tackled the case of the dimension $5$ (over a perfect field of characteristic other than $2$) in \cite{DU}. Nevertheless, some of their results are incorrect, as explained in \cite{MS}, and their classification lacks some Lie algebras. Maletesta and Siciliano solved the case of the dimension $5$ (over an algebraically closed field of characteristic $p>3$) in \cite{MS} by using a cohomological method that can be described as follows. As a first result, they showed that any $n$-dimensional $p$-nilpotent restricted Lie algebra can be obtained as a central extension by a restricted $2$-cocycle of a $(n-1)$-dimensional $p$-nilpotent restricted Lie algebra. Cohomologous cocycles lead to isomorphic Lie algebras. As a next step, they classified up to equivalence the restricted $2$-cocycles associated with the $p$-nilpotent restricted Lie algebras obtained in \cite{SU} and elaborated their central extensions. As a final step, isomorphic algebras need to be removed from the list.
Our paper builds on their work and uses our restricted cohomology for restricted Lie superalgebras to classify $p$-nilpotent restricted  Lie superalgebras of dimension $3$ and $4$. Our main tool is the concept of central extension, extensively studied in literature; see  for example \cite{KL, N, OR} in the case of  Lie algebras, see \cite{BGLL, IK,Ma, YL, SZ2} in the case of Lie superalgebras, see \cite{EF, EFP,EF2,EF3} in the case of restricted Lie algebras.

\subsection{Organization of the paper and main results} In Section \ref{basics}, we review basic concepts related to restricted Lie superalgebras. In Section \ref{restricted-cohomology}, we give an alternative definition of restricted $2$-cochains and $3$-cochains (Definition \ref{restricted-cochains}). The main novelty is that our definition involves odd elements in the definition of the ``restricted" term of the cochains. We then adapt the formulas for the differential maps and prove that they indeed define a (small) cohomology complex (Theorem \ref{dod0}), by following the proof of \cite[Theorem 3.13]{E}. This allows us to build restricted cohomology spaces of order $1$ and $2$. In addition, we investigate interpretations of the restricted cohomology in terms of derivations and extensions (Section \ref{sectionextensions}). The main novelty here concerns restricted central extensions. We introduce a $\Z_2$-graded sub-complex (see Eq. \eqref{subcomplex}) and show that the restricted central extensions are controlled by the cohomology induced by this sub-complex (Theorem \ref{thmext}). The rest of the paper is devoted to the classification of $p$-nilpotent restricted  Lie superalgebras of dimension $3$ and $4$ (over an algebraically closed field of characteristic $p>2$). First, we examine the case of the dimension $3$ in Section \ref{sectiondim3}. Its main result is Theorem \ref{classif3}, where our method is not cohomological. Section \ref{sectiondim4} is devoted to the case of the dimension $4$. We show that any $p$-nilpotent restricted Lie superalgebra of dimension $n$ can be obtained as a central extension by a restricted $2$-cocycle of a $p$-nilpotent restricted Lie superalgebra of dimension $n-1$ (Proposition \ref{propext1dim}). Then, we classify up to equivalence the non-trivial \textit{ordinary} $2$-cocycles associated to nilpotent restricted Lie superalgebras of dimension $3$ (Theorem \ref{cocycles}), build the central extensions (Section \ref{extensions}) and remove isomorphic Lie superalgebras from the list (Theorem \ref{classif4}). The last step of our procedure is to investigate the $p$-nilpotent $p|2p$ structures on these $4$-dimensional Lie superalgebras, which is achieved in Theorem \ref{pmap4} using Jacobson's Theorem \ref{SJac}. Finally, we compare our method (on an example) to the ``purely cohomological method" which consists in building extensions with \textit{restricted} cocycles and computing new $p|2p$-maps using those cocycles (Section \ref{comparison}). We conclude the paper by investigating $p|2p$-structures on some higher-dimensional nilpotent Lie superalgebras that were obtained in \cite{GKN} (Section \ref{knm}). \\

\noindent\textbf{Conventions and notations.} In the paper, $\mathbb{K}$ is an arbitrary field of characteristic $\mathrm{char}(\mathbb{K})=p>2$.  ``Ordinary" means ``not restricted" and the group of integers modulo $2$ is denoted by $\Z_{2}$. Let $V$ be a $\Z_{2}$-graded space. The degree of an homogeneous element $v\in V_{\bar{i}}$ is denoted by $|v|:=\bar{i}$. The element $v$ is called \textit{even} if $v\in V_\ev$ and \textit{odd} if $v\in V_\od$. 

\section{Restricted Lie (super)algebras}\label{basics}
For a comprehensive study of restricted Lie algebras, see \cite{SF}.

\subsection{Restricted Lie algebras}\label{defres}

Let $L$ be  a~finite-dimensional Lie algebra over~$\mathbb{K}$. Following \cite{J, SF}, a~map $(\cdot)^{[p]}:L\rightarrow L, \quad x\mapsto x^{[p]}$ is called a~\textit{$p$-structure} on $L$ and $L$ is said to be {\it restricted}  if 
\begin{equation}\label{RRR}
\begin{array}{l}
(i)\text{ $\ad_{x^{[p]}}=(\ad_x)^p$ for all $x\in L$;}\\[2pt]

(ii)\text{ $(\lambda x)^{[p]}=\lambda^p x^{[p]}$ for all $x\in L$ and for all $\lambda \in \mathbb{K}$;}\\[2pt]

(iii)\text{ $(x+y)^{[p]}=x^{[p]}+y^{[p]}+\displaystyle\sum_{1\leq i\leq p-1}s_i(x,y)$, where the  $s_i(x,y)$ can be obtained from}\\[2mm]
(\ad_{\lambda x+y})^{p-1}(x)=\displaystyle\sum_{1\leq i \leq p-1} is_i(x,y) \lambda^{i-1}.\\
\end{array}
\end{equation}

As a convenience, an explicit expression of  $\displaystyle\sum_{1\leq i\leq p-1}s_i(x,y)$ in terms of nested brackets is given by
\begin{align}\label{si}
\sum_{i=1}^{p-1}s_i(x,y)&=\sum_{\underset{x_{p-1}=y,~x_p=x}{x_k\in\{x,y\}}}\frac{1}{\sharp(x)}[x_1,[x_2,[...,[x_{p-1},x_p]...]]],
\end{align}
where $\sharp(x):=\text{card}\{k,~x_k=x,~k=1,\cdots,p\}$. Observe that $1\leq\sharp(x)\leq p-1$ in the expression above.\\

    
To investigate $p$-structures on a Lie algebra, it is useful to use the following theorem, due to Jacobson.
\sssbegin{Theorem}[\cite{J}] \label{Jac}
Let $(e_j)_{j\in J}$ be a~basis of $L$ such that there are $f_j\in L$ satisfying $(\ad_{e_j})^p=\ad_{f_j}$. Then, there exists exactly one $p$-map  $(\cdot)^{[p]}:L\rightarrow L$ such that 
\[
e_j^{[p]}=f_j \quad \text{ for all $j\in J$}.
\]
\end{Theorem}

Let $L$ be a restricted Lie algebra. A $L$-module $M$ is called \textit{restricted} if
\[
\underbrace{x\cdots x}_{p\text{~~times}}\cdot m =x^{[p]}\cdot m \quad \text{for all $x\in L$ and any $m\in M$.}
\]

\subsection{Restricted Lie superalgebras}

Let $L$ be a~finite-dimensional  Lie superalgebra defined over a field of characteristic $p>2$. For $p=3$, the Jacobi identity does not imply $[x,[x,x]]=0$, for all $x\in L_\od$, therefore we require it as part of the definition. 

Following \cite{Mi,P}, we say that $L$ has a~\textit{$p|2p$-structure}\footnote{In \cite{Mi,P}, it is called a $p$-superalgebra. The terminology $p|2p$ was introduced in \cite{BGL09}.} if $L_\ev$ is  a restricted Lie algebra and 
\begin{equation}
\label{RRRS}
\text{$\ad_{x^{[p]}}(y)=(\ad_x)^p(y)$  for all $x \in L_\ev$ and for all $y\in L$.}
\end{equation}
We set 
\[
(\cdot)^{[2p]}:L_\od \rightarrow L_\ev, \quad x\mapsto (x^2)^{[p]}  \text{~~for any $x\in L_\od$, where } x^2:=\frac{1}{2}[x,x].
\]
The pair $\bigl(L, (\cdot)^{[p|2p]}\bigl)$ is referred to as a~\textit{restricted} Lie superalgebra. In the case where the base field is of characteristic $p=2$, there are various types of restrictedness (see \cite{BGL09}), but we don't consider them in this work.

The following theorem is a~straightforward superization of Jacobson's theorem \ref{Jac}.
\sssbegin{Theorem}\label{SJac}
Let $(e_j)_{j\in J}$ be a~basis of $L_\ev$, and let the elements $f_j\in L_\ev$ be such that  ${(\ad_{e_j})^p=\ad_{f_j}}$. Then, there exists exactly one $p|2p$-mapping $(\cdot)^{[p|2p]}:L\rightarrow L$ such that 
\[
e_j^{[p]}=f_j \quad \text{ for all $j\in J$}.
\]
\end{Theorem}

Let $\bigl(L, (\cdot)^{[p|2p]_L}\bigl)$ and $\bigl(H, (\cdot)^{[p|2p]_H}\bigl)$ be two restricted Lie superalgebras. A morphism of Lie superalgebras $\varphi: L\rightarrow H$ is called \textit{restricted} if in addition to $\varphi\bigl([x,y]\bigl)=\bigl[\varphi(x),\varphi(y)\bigl]~\forall x,y\in L$, we also have 
$$\varphi\bigl(x^{[p]_L}\bigl)=\varphi(x)^{[p]_H},~\forall x\in L_\ev.  $$

A homogeneous ideal $I=I_\ev\oplus I_\od$ of $L$ is called a~\textit{$p$-ideal} if it is closed under the $p|2p$-map; namely 
\[
\text{$x^{[p]}\in I_\ev$ for all $x\in I_\ev$}. 
\]
As $x^{[2p]}=(x^2)^{[p]}$, for all $x\in L_\od$, this would imply that 
\[
x^{[2p]}\in I_\ev \text{ for all $x\in I_\od$}.
\]

The \textit{center} of a restricted Lie superalgebra $L$ is the $p$-ideal 
$$\fz(L)=\bigl\{ x\in L,~ [x,y]=0~\forall y\in L \bigl\}.   $$

Let $L$ be a restricted Lie superalgebra. An $L$-module $M$ is called \textit{restricted} if
\[
\underbrace{x\cdots x}_{p\text{~~times}} \cdot m =x^{[p]}\cdot m \quad \text{for all $x\in L_\ev$ and $m\in M$.}
\]
\subsection{$p$-nilpotent restricted Lie (super)algebras}
Let $L$ be a Lie superalgebra. We define a \textit{descending central sequence} by
    $$C^0(L)=L,~\text{ and }~ C^{k+1}(L)=\bigl[C^k(L),L\bigl].$$
The Lie superalgebra $L$ is called \textit{nilpotent} if there exists $k\geq 0$ such that $C^k(L)=0$. The smallest $k\geq 0$ satisfying $C^k(L)=0$ is called the \textit{nilindex} of $L$ and $L$ is said to be nilpotent of order $k$.
Suppose that $L$ is restricted. Then $L$ is called $p$\textit{-nilpotent} if there exists $n\geq0$ such that $x^{[p]^n}=0~\forall x\in L_\ev,$ where $x^{[p]^n}=(\cdots(x^{[p]})^{[p]}\cdots)^{[p]}$.   Any $p$-nilpotent restricted Lie superalgebra is nilpotent, using Engel's theorem.

\section{Restricted cohomology for restricted Lie superalgebras} \label{restricted-cohomology}

Evans and Fuchs constructed the restricted cohomology of restricted Lie algebras up to order $2$ in the most general case (\cite{E,EF}). In \cite{YCC} an attempt is made to generalize this cohomology to restricted Lie superalgebra. Their cohomology is not suitable for encoding some natural constructions, such as extensions of superalgebras, due to the absence of odd terms in their definition (see Section \ref{sectionextensions}). Our goal here is to fix this definition by taking into account odd elements and the missing sign.

\subsection{Chevalley-Eilenberg cohomology for Lie superalgebras}
Let $L=L_\ev\oplus L_\od$ be a restricted Lie superalgebra and let $M=M_\ev\oplus M_\od$ be a restricted module. Let us review the Chevalley-Eilenberg cohomology for Lie superalgebras, see \cite{BGL} for more details. 
For $n=0$, we put $C^0_{\text{CE}}(L; M):=M$. For $n>0$, the space of cochains $C^n_{\text{CE}}(L, M)=\wedge^n(L;M)$ is the space of multilinear super anti-symmetric maps $L^{\times n}\rightarrow M$. The differential is given by:
\[
\begin{array}{l}
 d_{\mathrm{CE}}^0(m)(x)=(-1)^{|m||x|}x\cdot m\text{~~for all $m\in M$ and for all $x\in L$};\\[2mm]
d_{\mathrm{CE}}^n(\varphi)(x_1,\dots, x_{n+1}) \\
= \mathop{\sum}\limits_{i<j}(-1)^{|x_j|(|x_{i+1}|+\dots+|x_{j-1}|)+j}\varphi(x_1, \dots, x_{i-1}, [x_i,x_j], x_{i+1} \dots, \hat  x_j, \dots, x_{n+1})\\
+\mathop{\sum}\limits_{j}(-1)^{|x_j|(|\varphi|+|x_{1}|+\dots+|x_{j-1}|)+j}x_j \cdot \varphi (x_1, \dots, \hat  x_j, \dots, x_{n+1})\\
\text{~~for any $\varphi \in C^{n}_{\text{CE}}(L; M)$ with $n>0$, and $x_1,\dots, x_{n+1}\in L$.}  \\ \text{ Here, the symbol $(~\hat{\cdot}~)$ means that the term is omitted.}
\end{array}
\]
The spaces $C^n_{\ce}(L;M)$ are $\Z_2$-graded. The degree of an homogeneous cochain $\varphi\in C^n_{\ce}(L;M)$ is defined by
$$|\varphi|:=|\varphi(x_1,\cdots, x_n)|-|x_1|-|x_2|-\cdots-|x_n|,~\text{ for homogeneous }x_i\in L.$$ A Chevalley-Eilenberg cochain will sometimes be called \textit{ordinary cochain}.

\subsection{Restricted cohomology for Lie superalgebras}\label{restcohomosuper}
The restricted cohomology for restricted Lie algebras has been developed by Evans and Fuchs (\cite{E,EF}). For an arbitrary Lie algebra $L$, they provided explicit descriptions of the restricted cochains spaces $C^n_{*}(L;M)$ for $n=0,1,2,3$ and of the differentials of order $0,1$ and $2$. A generalization of this work to restricted Lie superalgebras will follow. 

\subsubsection{Restricted cochains.}\label{restricted-cochains} For $n\in\{0,1\}$, we put $C^n_{*}(L;M)=C^n_{\ce}(L;M)$.  Let $\varphi\in C^2_{\ce}(L;M)$ and $\omega: L_\ev\longrightarrow M$ be a map.
We say that $\omega$ is \textit{$\varphi$-compatible}\footnote{In \cite{E}, this property is called \textit{$(*)$-property}. The terminology ``$\varphi$-compatible" has been introduced in \cite{EF3}.} if
	\begin{enumerate}
		\item[$(i)$] $\omega(\lambda x)=\lambda^p\omega(x),~\forall \lambda\in \K,~\forall x\in L_\ev;$
		\item[$(ii)$] 
        \begin{align*}
        \omega(x+y)&=\omega(x)+\omega(y)+\displaystyle\sum_{\underset{x_1=x,~x_2=y}{x_i\in\{x,y\}}}\frac{1}{\sharp(x)}\varphi\Bigl([[\cdots[x_1,x_2],x_3]\cdots,x_{p-1}],x_{p}\Bigl)\\
        &+\displaystyle\sum_{\underset{x_1=x,~x_2=y}{x_i\in\{x,y\}}}\frac{1}{\sharp(x)}\sum_{k=1}^{p-2}(-1)^kx_p\cdots x_{p-k+1}\varphi\Bigl([[\cdots[x_1,x_2],x_3]\cdots,x_{p-k-1}],x_{p-k}\Bigl),
        \end{align*}
	\end{enumerate}
with $x,y\in L_\ev $ and where $\sharp(x):=\text{card}\{k,~x_k=x,~k=1,\cdots,p\}$.

We set 
$$ C^2_*(L;M):=\left\lbrace (\varphi,\omega),~\varphi\in C^2_{\ce}(L;M),~\omega: L_\ev\rightarrow M \text{ is $\varphi$-compatible }\right\rbrace.     $$

Let $\alpha\in C^3_{\ce}(L;M)$ and $\beta: L\times L_\ev\rightarrow M$ be a map \footnote{In \cite{YCC}, the map $\beta$ was defined on $L_\ev\times L_\ev$ which is not suitable to capture central extensions. }. We say that $\beta$ is \textit{$\alpha$-compatible} if
	\begin{enumerate}
		\item[$(i)$] $\beta(x,y)$ is linear with respect to $x$;
		\item[$(ii)$]$\beta(x,\lambda y)=\lambda^p\beta(x,y)$;
		\item[$(iii)$] \begin{align*}\beta (x,y_1+y_2)&=\beta(x,y_1)+\beta(x,y_2)\\&-\displaystyle\sum_{\underset{h_1=y_1,~h_2=y_2}{h_i\in\{y_1,y_2\}}}\frac{1}{\sharp(y_1)}\sum_{j=0}^{p-2}(-1)^j\\&\times\sum_{k=1}^{j}\binom{j}{k}h_p\cdots h_{p-k-1}\alpha\Bigl( [\cdots[x,h_{p-k}],\cdots,h_{p-j+1}],[\cdots[h_1,h_2],\cdots,h_{p-j-1}],h_{p-j}\Bigl),\end{align*} 
	\end{enumerate}
with $\lambda\in \K$, $x \in L,y,y_1,y_2\in L_\ev$ and $\sharp(y_1)=\text{card}\{i,~h_i=y_1,~i=1,\cdots,p\}$. We set 
$$ C^3_*(L;M):=\left\lbrace (\alpha,\beta),~\alpha\in C^3_{\ce}(L;M),~\beta: L\times L_\ev \longrightarrow M \text{ is $\alpha$-compatible } \right\rbrace.$$ Note that the first component of $\beta$ is allowed to be odd.\\

\subsubsection{Technical difficulties.}\label{tech} The spaces  $C^2_*(L;M)$ and  $C^3_*(L;M)$ are actually not naturally $\Z_2$-graded because the second component of a cochain is not linear. Indeed, it is possible that $\varphi\in C^2_{\ce}(L;M)_\ev$, whereas $\omega(x)\in M_\od$ for all $x\in L_\ev$. Here is an example.

Consider the Lie superalgebra $L={\bf L_{2|2}^5}=\left\langle x_1,x_2|x_3,x_4,~[x_1,x_3]=x_4\right\rangle$ (see Theorem \ref{classif4}), with a  $p|2p$-structure given by $x_1^{[p]}=x_2,~x_2^{[p]}=0.$ Consider the even ordinary $2$-cocyle $\varphi=x_3\otimes(x_1^*\wedge x_4^*)\in Z^2_{\ce}(L;L)_\ev$. Let us define a map $\omega:L_\ev\rightarrow L_\od$ by $\omega(x_1)=\omega(x_2)=x_4$. This map is well defined and is $p$-semilinear (\textit{i.e,} $\omega(\lambda x+y)=\lambda^p\omega(x)+\omega(y),~\forall x,y\in L,~\forall \lambda\in \K$) on $L_\ev$. Moreover, we have $\ind^2(\varphi,\omega)=0$ (see Section \ref{op}). Thus, $(\varphi,\omega)\in Z^2_*(L;L)$.

For $(\varphi,\omega)\in C^2_*(L;M),$ we write \begin{equation}\label{lin}(\varphi,\omega)=(\varphi_\ev,\omega_\ev)+(\varphi_\od,\omega_\od), \text{ where } \text{Im}(\omega_{\bar{j}})\subseteq M_{\bar{j}}.\end{equation} Observe that also $(\varphi_{\bar{j}},\omega_{\bar{j}})\in C^2_*(L;M),$ thanks to the $\varphi$-compatibility.

\subsubsection{Differential operators and restricted cocycles}\label{op}
We define a map $d^0_*:C^0_*(L;M)\rightarrow C^1_*(L;M)$ by $d^0_*=d^0_{\ce}$.
An element $\varphi\in C^1_*(L;M)$ induces a map $\ind^1(\varphi): L_\ev\rightarrow M$ given by
$\ind^1(\varphi)(x)=\varphi\left(x^{[p]}\right)-x^{p-1}\varphi(x).$ It has been shown in \cite[Lemma 3.7]{E} that $\ind^1(\varphi)$ is $d^1_{\ce}\varphi$-compatible, therefore $d^1_*(\varphi):=\left(d^1_{\ce}\varphi,\ind^1(\varphi)\right)\in C^2_*(L;M)$. Moreover, it has been shown in \cite[Theorem 3.9]{E} that $d^1_*\circ d^0_*=0$. The proofs can be repeated without any changes since there are no odd elements involved in the definition of $\ind^1(\varphi)$.\\

Let us denote the spaces $Z^k_*(L;M)=\Ker\left(d^k_*\right)$ and $B^k_*(L;M)=\text{Im}\left(d^{k-1}_*\right)$, for $k\geq1$. The above discussion proves that the quotient space $H^1_*(L;M)=Z^1_*(L;M)/B^1_*(L;M)$ is well defined.

An element $(\varphi,\omega)\in C_*^2(L;M)$ induces a map\footnote{The sign $(-1)^{|\varphi||x|}$ was missed in \cite{YCC}.} $\ind^2(\varphi,\omega):L\times L_\ev\rightarrow M$ defined by
	$$\ind^2(\varphi,\omega)(x,y)=\varphi\left( x,y^{[p]}\right) -\displaystyle\sum_{i+j=p-1}(-1)^iy^i \varphi\Bigl( [[\cdots[x,\overset{j\text{ terms}}{\overbrace{y],\cdots],y}]},y\Bigl)+(-1)^{|\varphi||x|}x\omega(y),$$ for $\text{Im}(\omega)\subseteq M_{|\varphi|}$, and then extended using \eqref{lin}.	
The following result is a superized version of \cite[Lemma 3.10]{E}.
\sssbegin{Lemma}
Let $(\varphi,\omega)\in C_*^2(L;M)$. The map $\ind^2(\varphi,\omega)$ is $d^2_{\ce}\varphi$-compatible. 
\end{Lemma}

\begin{proof}
 First, we point out that \cite[Lemmas 3.7, 3.11, 3.12, 3.14]{E} are purely combinatorial, so they are also valid in the super setting. Moreover, \cite[Lemma 3.8]{E} is also valid if the first component is odd. These lemmas will be used freely in the proof. Let $(\varphi,\omega)\in C_*^2(L;M)$. From the definition of $\ind^2(\varphi,\omega)$, it is easy to see that conditions $(i)$ and $(ii)$ of the $d^2_{\ce}\varphi$-compatibility are satisfied. For further use, let us write for $u,v,w\in L$,
 \begin{align}\label{ce2}
    d^2_{\ce}\varphi(u,v,w)
    &=\varphi\bigl([u,v],w\bigl)-(-1)^{|v||w|}\varphi\bigl([u,w],v\bigl)-\varphi\bigl(u,[v,w]\bigl)\\
    &-(-1)^{|u||\varphi|}u\cdot\varphi(v,w)+(-1)^{|v|(|\varphi|+|u|)}v\cdot\varphi(u,w)\nonumber\\
    &-(-1)^{|w|(|\varphi|+|u|+|v|)}w\cdot\varphi(u,v).\nonumber
 \end{align}
\noindent\textbf{Notations.} For $1\leq j \leq p$, we denote $P_j:=\{p-j+1,\cdots, p-1,p\}$. Then, let $y_{t}\in L,~t\in P_j$ be a family indexed by $P_j$. Let $a:=\card(A)$ and $b:=\card(B)$ . We will use the notations $y_{t_A}:=y_{t_1}y_{t_2}\cdots y_{t_a}$ and $y_{t_B}:=y_{t_{a+1}}y_{t_{a+2}}\cdots y_{t_{a+b}}$. We will also use the notation $[x_1,\cdots, x_k]:=[[\cdots[x_1,x_2],x_3],\cdots,x_k]$, for all $x_1,\cdots, x_k\in L$.\\

Let $(\varphi,\omega)\in C_*^2(L;M)$, $x\in L$ and $y_1,y_2\in L_\ev$. 
\begin{align}
    \ind^2&(\varphi,\omega)(x,y_1+y_2)\nonumber=\\
    &~\varphi\left(x,y_1^{[p]}\right)+\varphi\left(x,y_2^{[p]}\right)\label{a}\\
    &+\displaystyle\sum_{\underset{k_1=1,~k_2=2}{k_s\in\{1,2\}}}\frac{1}{\sharp(y_1)}\varphi\bigl(x,[y_{k_1},y_{k_2},\cdots,y_{k_p}]\bigl)\label{b}\\
    &-\sum_{i+j=p-1}(-1)^iy_1^i\varphi\Bigl([x,\underset{j\text{~terms~}}{\underbrace{y_1,\cdots, y_1}],y_1} \Bigl)\label{c}\\
    &-\sum_{i+j=p-1}(-1)^iy_2^i\varphi\Bigl([x,\underset{j\text{~terms~}}{\underbrace{y_2,\cdots, y_2}],y_2} \Bigl)\label{d}\\
    &-\sum_{i=0}^{p-1}\sum_{\underset{\sharp(y_1)<p,~\sharp(y_2)<p}{k_s\in\{1,2\}}}(-1)^iy_{k_1}\cdots y_{k_i}\varphi\bigl([x,y_{k_{i+1}},\cdots,y_{k_{p-1}}],y_{k_p}\bigl)\label{e}\\
    &+(-1)^{|\varphi||x|}\biggl(x\omega(y_1)+x\omega(y_2)\label{f}\\
    &~~~~~~~~+\displaystyle\sum_{\underset{k_1=1,~k_2=2}{k_s\in\{1,2\}}}\frac{1}{\sharp(y_1)}\displaystyle\sum_{j=0}^{p-2}x\cdot y_{k_p}\cdots y_{k_{p-j-1}}\varphi\bigl([y_{k_1},\cdots,y_{k_{p-j+1}}],y_{k_{p-j}}\bigl)\biggl).\label{g}
\end{align}
First, we have $\eqref{a}+\eqref{c}+\eqref{d}+\eqref{f}=\ind^2(\varphi,\omega)(x,y_1)+\ind^2(\varphi,\omega)(x,y_2)$. It remains to show that $\eqref{b}+\eqref{e}+\eqref{g}$ gives the double sum in Condition $(iii)$ of the $d^2_{\ce}\varphi$-compatibility. First, in the case where $|x|=\ev$, the proof of \cite[Lemma 3.10]{E} can be repeated word for word. From now on, we assume  that $|x|=\od$.
We focus on \eqref{g}. Every term is of the form $u\cdot\varphi\bigl(v,w\bigl)$, with $u=x\cdot y_{k_p}\cdots y_{k_{p-j-1}}$ being odd, $v=[y_1,\cdots,y_{k_{p-j+1}}]$ and $w=y_{k_{p-j}}$ being even.
Rewriting \eqref{ce2} with $|v|=|w|=\ev$, we get

\begin{align}\label{ce2bis}
u\cdot\varphi(v,w)&=(-1)^{|u||\varphi|}\biggl(\varphi\bigl([u,v],w\bigl)+\varphi\bigl([v,w],u\bigl)+\varphi\bigl([w,u],v\bigl)\\
&-v\cdot\varphi(w,u)-w\cdot\varphi(u,v)-d^2_{\ce}\varphi(u,v,w)\biggl).\nonumber
\end{align}
A careful computation shows that
    \begin{align*}
        \displaystyle\sum_{u,v,w}d^2_{\ce}\varphi(u,v,w)
        &:=\displaystyle\sum_{A\sqcup B=P_j}d^2_{\ce}\varphi(x\cdot y_{k_B},y_{k_A},y_{k_p-j})\\
        &=\displaystyle\sum_{\underset{l_1=1,~l_2=2}{l_s\in\{1,2\}}}\frac{1}{\sharp(y_1)}\sum_{j=0}^{p-2}(-1)^j \sum_{k=1}^{j}\binom{j}{k}y_{l_{p}}\cdots y_{l_{p-k-1}}\\
        & ~\hfill~\times d^2_{\ce}\bigl([x,y_{l_{p-k}},\cdots,y_{l_{p-j+1}}],[y_{l_{1}},\cdots y_{l_{p-j-1}}],y_{l_{p-j}}\bigl),
    \end{align*}
which is exactly the term that appears in Condition $(iii)$ of the $d^2_{\ce}\varphi$-compatibility, up to a sign. The sign of the expression $d^2_{\ce}\varphi(u,v,w)$ in \eqref{ce2bis} is $-(-1)^{|x||\varphi|}$. By plugging it into \eqref{g}, we have to multiply again by $(-1)^{|x||\varphi|}$, therefore the sign is just $(-1)$ and we get precisely Condition $(iii)$.

For the remaining terms, namely \eqref{b}, \eqref{e} and part of \eqref{g}, the same argument apply and a cautious analysis shows that all the signs $(-1)^{|x||\varphi|}$ cancel each other. Therefore, we can apply the same arguments as in Evans' proof, namely we obtain \cite[Eq. (3.34)-(3.40)]{E}. The following case-by-case study remain valid. The end of Evans' proof is therefore valid and the remaining terms cancel. The proof is complete.
\end{proof}

Let us define the differential map $d^2_*:C^2_*(L;M)\rightarrow C^3_*(L;M)$ by $d^2_*(\varphi,\omega)=\left(d^2_{\ce}\varphi,\ind^2(\varphi,\omega)\right)$.

\sssbegin{Theorem}\label{dod0}
    We have $d^2_*\circ d^1_*=0$. Therefore, the quotient space $$H^2_*(L;M)=Z^2_*(L;M)/B^2_*(L;M)$$ is well defined.
\end{Theorem}

\begin{proof}
For further use, let us write
\begin{equation}\label{ce1}
        d^1_{\ce}\psi(u,v)=  \psi\bigl([u,v]\bigl)-(-1)^{|u||\psi|}u\cdot\psi(v)+(-1)^{|v|(|\psi|+|u|)}v\cdot\psi(u),~\forall u,v\in L.
\end{equation}
It is well-known that $d^2_{\ce}\circ d^1_{\ce}=0$.  Let $\psi\in C^1_*(L;M)$. It remains to show that $\ind^2\left(d^1_{\ce}(\psi),\ind^1(\psi)\right)=0.$ Let  $x,y\in L$. In the case where both $x$ and $y$ are even, the proof of \cite[Theorem 3.13]{E} can be repeated word for word. Let us suppose that $x\in L_\od$ and that $y\in L_\ev$.
\begin{align}
    \ind^2&\left(d^1_{\ce}(\psi),\ind^1(\psi)\right)(x,y)=\nonumber\\
        &~\psi\bigl(x,y^{[p]}\bigl)\underline{-(-1)^{|\psi|}x\psi\left(y^{[p]}\right)}+y^{[p]}\psi(x)\label{aa}\\
        &-\sum_{i+j=p-1}(-1)^iy^i\biggl( \psi\bigl( [x,\overset{j\text{ terms}}{\overbrace{y,...,y}]},y\bigl)\label{bb}\\
        &~~~~~~-(-1)^{|\psi|}[x,\overset{j}{\overbrace{y,\cdots,y}]}\psi(y)\label{cc}\\ &~~~~~~+y\psi\bigl([x,\overset{j}{\overbrace{y,\cdots,y}]}]\bigl)\biggl)\label{dd}\\
        &+\underline{(-1)^{|\psi|}\Bigl(x\psi\left(y^{[p]}\right)}-xy^{p-1}\psi(y)\Bigl)\label{ee}.
\end{align}

The underlined terms cancel. The term $y^{[p]}\psi(x)$ of \eqref{aa} cancels with \eqref{dd} when $j=0$, because $M$ is a restricted module. A similar argument shows that the term $\psi\bigl(x,y^{[p]}\bigl)$ of \eqref{aa} cancels with \eqref{bb} if $j=p-1$. The other terms of \eqref{bb} and \eqref{dd} cancel each other. The last remaining terms are \eqref{cc} and \eqref{ee}, that can be factorized by $(-1)^{|\psi|}$. We finish the proof by repeating the arguments of \cite[Theorem 3.13]{E}.
\end{proof}

We are now in the following situation.
$$0\longrightarrow C^0_*(L;M)\overset{d^0_*}{\longrightarrow}C^1_*(L;M)\overset{d^1_*}{\longrightarrow}C^2_*(L;M)\overset{d^2_*}{\longrightarrow} C^3_*(L;M).$$
We expect that extending the complex to higher orders and providing explicit differential formulas will be very challenging.

\subsubsection{Example} Here we provide an example of computation of restricted cohomology with adjoint coefficients. Consider the Lie superalgebra $L=\bf{L_{1|2}^3}$ of superdimension $(1|2)$ given in the basis $e_1|e_2,e_3$ (even $|$ odd) by the bracket $[e_1,e_2]=e_3$ and the $p|2p$-map $e_1^{[p]}=0.$ Let $\varphi\in C^2_{\ce}(L;L)$. In this case, a map $\omega:L_\ev\rightarrow L$ is $\varphi$-compatible if and only if
\begin{equation}\label{psemi}\omega(\lambda x)=\lambda^p\omega(x) \text{ and } \omega(x+y)=\omega(x)+\omega(y),~\forall x,y\in L_\ev,~\forall \lambda\in \K.  \end{equation}
We will use the notation $\Delta_{i,j}$ to write down the Chevalley-Eilenberg $2$-cocycles, see Eq. \eqref{deltacocycles}.

 A basis for the Chevalley-Eilenberg $2$-cocycles space $Z_{\ce}^2(L;L)$ is given by

\[
\begin{array}{rlllll}
\varphi_1 &= &e_1\otimes \Delta_{1,2}+2e_2\otimes\Delta_{2,2}; &\varphi_2 &= &-2e_1\otimes \Delta_{1,3}+2e_3\otimes\Delta_{3,3};  \\[2mm] 
\varphi_3 &= &e_2\otimes\Delta_{1,2}; &\varphi_4 &= &e_2\otimes\Delta_{1,3};  \\[2mm]
\varphi_5 &= &2e_2\otimes\Delta_{2,2}+e_3\otimes\Delta_{2,3}; &\varphi_6 &= &e_3\otimes\Delta_{1,2};  \\[2mm]
\varphi_7 &= &e_3\otimes\Delta_{1,3}; &\varphi_8 &= &e_3\otimes\Delta_{2,2}.  \\[2mm]
\end{array}
\]
\underline{The case where $p>3$.} Let $\varphi\in Z_{\ce}^2(L;L)$, and let $\omega$ be a map satisfying \eqref{psemi}. Then, $(\varphi,\omega)\in Z_{*}^2(L;L)$ if and only if $\omega(e_1)=\gamma e_3$ for some $\gamma\in\K$. Moreover, one can check that $\varphi_6$ and $\varphi_8$ are Chevalley-Eilenberg coboundaries, that $\varphi_5\cong\varphi_1, \varphi_7\cong\varphi_3$  and that $\ind^1(\psi)=0~\forall \psi\in C^1_*(L;L)$. Therefore, we have
$$H^2_*(L;L)=\Span\bigl((\varphi_1,0);~(\varphi_2,0);~(\varphi_3,0);~(\varphi_4,0);~(0,\omega_5)\bigl),   $$ where $\omega_5(e_1)=e_3$.

\underline{The case where $p=3$.} Let $\omega$ be a map satisfying \eqref{psemi}, given by $\omega(e_1)=\gamma_1e_1+\gamma_2e_2+\gamma_3e_3,~\gamma_1,\gamma_2,\gamma_3\in\K$. For every $\varphi\in Z^2_{\ce}(L;L)$, we have $\ind^2(\varphi,\omega)(e_1,e_1)=\gamma_2 e_3.$ Therefore, $\gamma_2=0.$ For $i\neq 4$, we have $\ind^2(\varphi_i,\omega)(e_2,e_1)=\gamma_1 e_3$.  Moreover, we have 
$$\ind^2(\varphi_4,\omega)(e_2,e_1)=[\varphi_4(e_3,e_1),e_1]-\gamma_1e_3=(1-\gamma_1)e_3. $$
Therefore, $(\varphi_4,\omega)\in Z^2_*(L;L)$ if and only if $\gamma_1=1$. As a consequence,
$$H^2_*(L;L)=\Span\bigl((\varphi_1,0);~(\varphi_2,0);~(\varphi_3,0);~(\varphi_4,\omega_4);~(0,\omega_5)\bigl),   $$ where 
$\omega_4(e_1)=e_1$ and $\omega_5(e_1)=e_3$.\\

The results of Sections \ref{outerderi} and \ref{modulesextensions} below were already proven in \cite{YCC}, but we have included them here for completeness.

\subsection{Restricted outer derivations}\label{outerderi} Let $L$ be a restricted Lie superalgebra. A \textit{restricted derivation} of $L$ is an homogeneous map $D:L\rightarrow L$ satisfying
\begin{enumerate}
    \item[(i)] $D\bigl([x,y]\bigl)=[D(x),y]+(-1)^{|x||D|}[x,D(y)],~\forall x,y\in L;$
    \item[(ii)] $D\bigl(x^{[p]}\bigl)=\ad_x^{p-1}\circ D(x),~\forall x\in L_\ev.$
\end{enumerate}
We denote $\Der_p(L)$ the space of restricted derivations of the restricted Lie superalgebra $L$. This space has a natural $\Z_2$-graduation given by the degree of the elements. The subspace $\ad(L)=\{\ad_x,~x\in L\}$ is a restricted ideal of $\Der_p(L)$.

\sssbegin{Theorem}\textup{(}\cite[Theorem 3.2]{YCC}\textup{)}
    Let $L$ be a restricted Lie superalgebra. Then $\Der_p(L)/\ad(L)=H^1_*(L;L)$.
\end{Theorem}
\begin{proof}
See \cite[Theorem 3.2]{YCC}.
\end{proof}
\subsection{Restricted extensions of restricted modules}\label{modulesextensions} Let $L$ be a restricted Lie superalgebra and let $M,N$ be two restricted $L$-modules. A restricted extension of $N$ by $M$ is a short exact sequence of restricted $L$-modules
\begin{equation}\label{extmod}
    0\longrightarrow M\overset{\iota}{\longrightarrow}E\overset{\pi}{\longrightarrow}N\longrightarrow 0.
\end{equation}
Two extensions $E_1$ and $E_2$ of $N$ by $M$ are called equivalent if there is a restricted isomorphism of $L$-modules $\sigma:E_1\rightarrow E_2$ such that the usual diagram commutes, see \cite[Definition 3.3]{YCC}. The space $\Hom(N,M)$ of restricted $L$-modules morphisms is a restricted $L$-module itself with the action
$$(x\cdot\phi)(n):=x\cdot\phi(n)+(-1)^{|\phi||x|}\phi(x\cdot n),~\forall x\in L,~ \forall\phi\in\Hom(N,M),~\forall n\in N.  $$

\sssbegin{Theorem}\textup{(}\cite[Theorem 3.4]{YCC}\textup{)}
 Let $L$ be a restricted Lie superalgebra and let $M,N$ be two restricted $L$-modules. The equivalence classes of restricted extensions of $N$ by $M$ is in one to one correspondence with $H^1_*\bigl(L;\Hom(N,M)\bigl)$.
\end{Theorem}
\begin{proof}
See \cite[Theorem 3.4]{YCC}.
\end{proof}

\subsection{Restricted extensions of restricted Lie superalgebras}\label{sectionextensions}

In \cite{YCC}, restricted central extensions are discussed. However, the main result (\cite[Theorem 3.7]{YCC}) is not correct. Here is an attempt to rectify it.
Let $\bigl(L,[\cdot,\cdot],(\cdot)^{[p]}\bigl)$ be a restricted Lie superalgebra, and $M$ be a strongly abelian restricted Lie superalgebra (\textit{i.e}, $[m,n]=0~\forall m,n\in M,~\text{and }m^{[p]}=0~\forall m\in M_\ev$). A \textit{restricted extension} of $L$ by $M$ is a short exact sequence of restricted Lie superalgebras
\begin{equation}\label{ext}
    0\longrightarrow M\overset{\iota}{\longrightarrow}E\overset{\pi}{\longrightarrow}L\longrightarrow 0.
\end{equation}
Note that in the above diagram, the maps $\iota$ and $\pi$ are \textit{even}.
Two restricted central extensions of $L$ by $M$ are called \textit{equivalent} if there is a restricted Lie superalgebras morphism $\sigma:E_1\rightarrow E_2$ such that the following diagram commutes:
\begin{equation}\label{eqext}
\begin{tikzcd}
            &                                                & E_1 \arrow[rd, "\pi_1"] \arrow[dd, "\sigma"] &             &   \\
0 \arrow[r] & M \arrow[ru, "\iota_1"] \arrow[rd, "\iota_2"'] &                                              & L \arrow[r] & 0. \\
            &                                                & E_2 \arrow[ru, "\pi_2"']                     &             &  
\end{tikzcd}
\end{equation}

\subsubsection{Technicalities.} Let $L$ be a restricted Lie superalgebra and consider a restricted extension of $L$ given by a short exact sequence \eqref{ext}. Then, $\iota(M)$ is an (homogeneous) $p$-ideal of $E$ and $M$ is a restricted $L$-module. Indeed, let $m\in M$ and $x\in E$. Then $[\iota(m),x]\in\Ker(\pi)=\text{Im}(\iota),~\forall x\in E$. Now, let $c=c_\ev+c_\od\in \iota(M)$, with $c_\ev\in E_\ev$ and $c_\od\in E_\od$. We have $\pi(c_\ev)=\underset{=0}{\underbrace{\pi(c_\ev+c_\od)}}-\pi(c_\od)$. We deduce that $\pi(c_\ev)=\pi(c_\od)=0$, therefore $c_\ev\in \iota(M)_\ev$, $c_\od\in \iota(M)_\od$ and $\iota(M)$ is homogeneous.\\

We define a restricted $L$-module structure on $M$ with
\begin{equation}\label{action2}
    x\cdot m:=\iota^{-1}\Bigl( [\Tilde{x},\iota(m)]\Bigl)\in M,~\forall x\in L,~\forall m\in M.
\end{equation}
Let $x\in L_\ev$. It exists $\tilde{x}\in E_\ev$ such that $\pi(\tilde{x})=x$ (\textit{a priori}, there is no guarantee that $\tilde{x}$ is even). Suppose that $\tilde{x}=\tilde{x}_\ev+\tilde{x}_\od$, with $\tilde{x}_\ev\in E_\ev$ and $\tilde{x}_\od \in E_\od$. Then, $\tilde{x}_\od \in \Ker(\pi)$. Thus, we can always choose $\Tilde{x}$ to be even, provided that $x$ is even. Therefore, we can consider elements of the form $\Tilde{x}^{[p]}$, which allows us to prove that the $L$-module $M$ is indeed restricted. Let $x\in L_\ev$, $m\in M$, and choose $\Tilde{x}\in E_\ev$ such that $\pi(\Tilde{x})=x$.

First, we have $\pi\bigl(\Tilde{x}^{[p]}\bigl)=\pi(\Tilde{x})^{[p]}=x^{[p]}=\pi\bigl(\widetilde{x^{[p]}}\bigl)$. Therefore, it exists $t_x\in\Ker(\pi)$ such that $\Tilde{x}^{[p]}=\widetilde{x^{[p]}}+t_x$. Since $\Ker(\pi)=\text{Im}(\iota)$, it exists $m_x\in M$ such that $t_x=\iota(m_x)$. Thus,
\begin{equation}\label{eqqqqq}
[t_x,\iota(m)]=[\iota(m_x),\iota(m)]=\iota\bigl([m_x,m]\bigl)=0, \text{ since $M$ is abelian.}
\end{equation}
Therefore, using Equation \eqref{eqqqqq}, we have
\begin{align*}
    x^{[p]}\cdot m &=\iota^{-1}\Bigl(\bigl[\widetilde{x^{[p]}},\iota(m)\bigl]\Bigl) =\iota^{-1}\Bigl(\bigl[\tilde{x}^{[p]},\iota(m)\bigl]\Bigl) =\iota^{-1}\Bigl( \ad_{\Tilde{x}}^p\bigl(\iota(m)\bigl)\Bigl)\\
                   &=\iota^{-1}\Bigl(\bigl[\underset{p \text{ terms}}{\underbrace{\Tilde{x},[\cdots,[\Tilde{x}}},\iota(m)]\cdots\bigl] \Bigl) =\iota^{-1}\Bigl(\bigl[\Tilde{x},\iota\circ\iota^{-1}\bigl(\underset{p-1 \text{ terms}}{\underbrace{\Tilde{x},[\cdots,[\Tilde{x}}},\iota(m)]\bigl)\cdots\bigl] \Bigl)\\
                   &=x\cdot\iota^{-1}\Bigl(\ad^{p-1}_{\Tilde{x}}\circ\iota(m)\Bigl).
\end{align*}
By induction, we obtain $x^{[p]}\cdot m=\underset{p}{\underbrace{x\cdot x\cdots x}}\cdot m$, so the $L$-module $M$ is indeed restricted.

In the case where $\iota(M)\subset \fz(E)$, the action \eqref{action2} is identically zero and $M$ becomes a trivial $L$-module. These extensions are called \textit{restricted central extensions}.\\

Let $L$ be a restricted Lie superalgebra and $M$ a restricted $L$-module. We define a subspace  $C^2_*(L;M)^+\subset C^2_*(L;M)$ by
\begin{equation}\label{subcomplex}
    C^2_*(L;M)^+:=\Bigl\lbrace (\alpha,\beta)\in C^2_*(L;M),~\text{Im}(\beta)\subseteq  M_\ev  \Bigl\rbrace.
\end{equation}
The inclusion $C^2_*(L;M)^+\subset C^2_*(L;M)$ is strict in general, see Example \ref{tech}. Let us recall that $C^1_*(L;M)=C^1_{\ce}(L;M)$. In particular, the space $C^1_*(L;M)$ is $\Z_2$-graded and we have $(d^1_{\ce}\psi,\ind^1(\psi))\in C^2_*(L;M)^+$ if and only if $\psi\in C^1_*(L;M)_\ev$ or $\ind^1(\psi)=0$. 

\sssbegin{Lemma}\label{shapecocycles}
    Let $L$ be a restricted Lie superalgebra and $M$ a restricted $L$-module.
        \begin{itemize}
            \item[(i)] We have an inclusion $B^2_*(L;M)_\ev:=d^1_*\bigl( C^1_*(L;M)_\ev\bigl)\subset C^2_*(L;M)^+$.
            \item[(ii)] The space $C^2_*(L;M)^+$ is $\Z_2$-graded and the degree of an homogeneous element $(\alpha,\beta)\in C^2_*(L;M)^+$ is given by $|(\alpha,\beta)|=|\alpha|$ and $|(0,\beta)|=\ev$.
        \end{itemize}
\end{Lemma}
This Lemma allows us to consider the space $Z^2_*(L;M)^+:=\Ker\bigl(d^2_{*|C^2_*(L;M)^+}\bigl)$. Thus we can define 
$$H^2_*(L;M)^+:=Z^2_*(L;M)^+/B^2_*(L;M)_\ev.$$ The space $H^2_*(L;M)^+$ is $\Z_2$-graded.

We will show that $H^2_*(L;M)^+_\ev$ captures all restricted central extensions of the restricted Lie superalgebras by strongly abelian superalgebras. Compare this to \cite[Theorem 3.7]{YCC}, where their cohomology groups are actually not in one-to-one correspondence with restricted central extensions.

\sssbegin{Theorem}\label{thmext}
Let $L$ be a restricted Lie superalgebra and $M$ a strongly abelian restricted Lie superalgebra. Then, the equivalence classes of restricted central extensions of $L$ by $M$ are classified by $H^2_*(L;M)^+_\ev$. 
\end{Theorem}

\begin{proof} The proof consists of four steps.\\
   \noindent\underline{Step $1$}. Consider a central extension of $L$ by $M$ as in \eqref{ext}. It exists an even linear section $\rho:L\rightarrow E$ such that $\pi\circ\rho=\id_L$. One can show that
   $$\bigl[\rho(x),\rho(y)\bigl]-\rho\bigl([x,y]\bigl)\in\Ker(\pi) \text{ and } \rho(z)^{[p]}-\rho\bigl( z^{[p]}  \bigl)\in\Ker(\pi),~\forall x,y\in L,~z\in L_\ev.$$
    We define two maps
    \begin{align}
        \varphi: L\times L&\rightarrow M,~ (x,y)\mapsto \iota^{-1}\Bigl([\rho(x),\rho(y)]-\rho([x,y])\Bigl);\\
        \omega: L_\ev&\rightarrow M_\ev,~x\mapsto\iota^{-1}\Bigl(\rho(x)^{[p]}-\rho\bigl( x^{[p]}  \bigl)\Bigl).
    \end{align}
  Then, $(\varphi,\omega)\in Z^2_*(L;M)^+_\ev$. Indeed, the map $\omega$ is $\varphi$-compatible (see \cite[Lemma 3.25]{E}). Moreover, we have $d^2_{\ce}\varphi=0$ (see \cite[Theorem 3.7, (1.1)]{YCC}). It remains to show that $\ind^2(\varphi,\omega)(x,y)=0,~\forall x\in L,~\forall y\in L_{\ev}.$\footnote{In \cite[Theorem 3.7]{YCC}, both $x$ and $y$ are even, which is not correct.} Because $M$ is a trivial module, it is enough to show that 
\begin{equation}\label{cocococo}
\varphi\bigl(x,y^{[p]}\bigl)=\varphi\bigl([x,\underset{p-1}{\underbrace{y,\cdots,y}}],y\bigl)~\forall x\in L,~\forall y\in L_{\ev}.
\end{equation}
The map $\rho$ being even, the parity of the element $x$ is not involved in the computation. Thus, the proof of \eqref{cocococo} is the same as in \cite[Theorem 3.7, (1.2)]{YCC}. Therefore, $d^2_*(\varphi,\omega)=0$ and $(\varphi,\omega)\in Z^2_*(L;M)^+_\ev$.\\

\noindent\underline{Step $2$}. Let $(\varphi,\omega)\in Z_*^2(L;M)^+_\ev$. We define the two following maps on $L\oplus M$.
\begin{align}
   [x+m,y+n]_N&:=[x,y]+\varphi(x,y),~\forall x,y\in L,~ \forall m,n\in M;\label{extbracket}\\
    (x+m)^{[p]_N}&:=(x)^{[p]}+\omega(x),~\forall x\in L_\ev,~\forall m\in M_\ev.\label{extp}
\end{align} The next step consists in showing that the maps \eqref{extbracket} and \eqref{extp} define a $p|2p$-structure on $Z_*^2(L;M)^+_\ev$.
 To check that the map \eqref{extbracket} defines a Lie superalgebra structure is routine. 
    Since $\omega$ is $\varphi$-compatible, it is easy to see that $[\cdot,\cdot]_N$ and $(\cdot)^{[p]_N}$ satisfy  the points $(ii)$ and $(iii)$ of Definition \ref{defres}. It remains to check Condition \eqref{RRRS}. Let $x+m$ be a homogeneous element of $L\oplus M$ and $y+n\in (L\oplus M)_\ev$. Condition \eqref{RRRS} reads
    \begin{equation}\label{eqext2}
        \Bigl[x+m,(y+n)^{[p]_N}\Bigl]_N=\bigl[x,\underset{p}{\underbrace{y+n,\cdots,y+n}}\bigl]_N.
    \end{equation}
    The left-hand side of Eq. \eqref{eqext2} gives
    $$\bigl[x+m,(y+n)^{[p]_N}\bigl]_N=\bigl[x,y^{[p]}\bigl]+\varphi\bigl(x,y^{[p]}\bigl).$$
    The right-hand side of Eq. \eqref{eqext2} gives
    $$\bigl[x,\underset{p}{\underbrace{y+n,\cdots,y+n}}\bigl]_N=\bigl[x,\underset{p}{\underbrace{y,\cdots,y}}\bigl]+\varphi\bigl([x,\underset{p-1}{\underbrace{y,\cdots,y}}],y  \bigl).$$
Since $(\cdot)^{[p]}$ is a $p$-map, we have $\bigl[x,y^{[p]}\bigl]=\bigl[x,\underset{p}{\underbrace{y,\cdots,y}}\bigl]$ , hence \eqref{eqext2} is equivalent to
\begin{equation}
      \varphi\bigl(x,y^{[p]}\bigl)=\varphi\bigl([x,\underset{p-1}{\underbrace{y,\cdots,y}}],y \bigl),
\end{equation}
which is equivalent to $\ind^2(\varphi,\omega)=0$ provided that $M$ is a trivial module. \footnote{In the proof of \cite[Theorem 3.7, (2.2)]{YCC}, the authors consider a restricted cochain $(\alpha,\beta)$ to build the central extension. They claim that since $\ind^2(\alpha,\beta)=0$, the extension is equipped with a $p|2p$-structure. This is unfortunately not correct. Indeed, their definition of the map $\ind^2$ does not involve odd elements, making Condition $\eqref{RRRS}$ impossible to check. In the case of purely even modules, however, their result remains valid. In that case, the expressions of the form $``\alpha(\text{odd},\text{even})"$ vanish, so the problem does not occur.}\\

\noindent\underline{Step $3$}. Let $E_1$ and $E_2$ be two equivalent restricted central extensions of $L$ by $M$, built with the cocycles $(\varphi_1,\omega_1)$ and $(\varphi_2,\omega_2)$ respectively. We denote by $\sigma:E_1\rightarrow E_2$ the isomorphism of restricted Lie superalgebras that realizes the equivalence (see \eqref{eqext}). We will show that $(\varphi_1,\omega_1)$ and $(\varphi_2,\omega_2)$ are cohomologous. Denote $\rho_1$ (resp. $\rho_2$) the section of $\pi_1$ (resp. $\pi_2$). We have $\pi(\sigma\rho_1-\rho_2)=0$, therefore we can define an even map $\psi:L\rightarrow M,~x\mapsto\iota_2^{-1}(\sigma\rho_1-\rho_2)$.
As in the ordinary case, we have $\varphi_1=\varphi_2+d^1_{\ce}\psi$. Moreover, as in \cite[Lemma 3.25]{E}, we have $\ind^1(\psi)=\omega_2-\omega_1$. Hence $(\varphi_2-\varphi_1,\omega_2-\omega_1)\in B^2_*(L;M)_\ev$.\\

\noindent\underline{Step $4$}. Let $\psi: L\rightarrow M$ be an even linear map and $(\varphi_1,\omega_1)\in Z^2_*(L;M)^+_\ev$. Define $(\varphi_2,\omega_2)=(\varphi_1,\omega_1)+d^1_*(\psi)$. This fourth and final step of the proof consists in showing that the restricted central extensions $E_1$ and $E_2$, given by $(\varphi_1,\omega_1)$ resp. $(\varphi_2,\omega_2)$ are equivalent. We consider the bijective linear map $\sigma: E_1\rightarrow E_2$ defined by $\sigma(x+m)=x+m+\psi(x).$ This map makes the diagram \eqref{eqext} commute. It remains to show that $\sigma$ is a restricted morphism of restricted Lie superalgebras. Let $x_1+m_1$ and $x_2+m_2$ be homogeneous elements of $E_2$. Then,

\begin{align*}
    \bigl[\sigma(x_1+m_1),\sigma(x_2+m_2)\bigl]&=\bigl[x_1+m_1+\psi(x_1),x_2+m_2+\psi(x_2)\bigl]\\
        &=[x_1,x_2]+\varphi_2(x_1,x_2)\\
        &=[x_1,x_2]+\varphi_1(x_1,x_2)+\psi\bigl([x_1,x_2]\bigl)\\
        &=\sigma\bigl([x_1,x_2]+\varphi_1(x_1,x_2)   \bigl)\\
        &=\sigma\bigl([x_1+m_1,x_2+m_2]  \bigl).
\end{align*}
Let $x+m$ be an even element of $E_1$. Then,
$$\sigma\bigl((x+m)^{[p]}\bigl)=\sigma\bigl(x^{[p]}+\omega_1(x)\bigl)=x^{[p]}+\underset{=\omega_2(x)}{\underbrace{\omega_1(x)+\psi\bigl(x^{[p]}\bigl)}}=\bigl(x+m+\omega_2(x)\bigl)^{[p]}=\sigma(x+m)^{[p]}. $$
We conclude that $\sigma$ is a restricted isomorphism of restricted Lie superalgebras and the extensions $E_1$ and $E_2$ are therefore equivalent.
\end{proof}

\sssbegin{Corollary}
Let $L$ be a restricted Lie superalgebra. Then, the equivalences classes of one-dimensional restricted central extensions of $L$ by $M$ are classified by $H^2_*(L;\K)^+_\ev$.
\end{Corollary}

\section{$p$-nilpotent restricted  Lie superalgebras of dimension $3$}\label{sectiondim3}
    This section discusses non-isomorphic $p$-nilpotent restricted Lie superalgebras of  dimension $3$. Here is how we approached the problem. Let $L$ be a nilpotent Lie superalgebra of superdimension $\sdim(L)=(n|k)$ where $n+k=3$. The even part $L_\ev$ is a nilpotent Lie algebra of dimension $n$. Starting from each nilpotent Lie algebra of dimension $n$, we add a  new basis element and compute a new multiplication table. The list is then reduced by removing isomorphic Lie superalgebras. The majority of computations were performed with Mathematica, a computer algebra system.

\ssbegin{Lemma}
   Let $L$ be a nilpotent Lie superalgebra of total dimension $3$ over an
algebraically closed field $\K$ of characteristic $p\geq3$ . Then, $L$ is isomorphic to one of the following Lie superalgebras:\\

 \item \underline{$\sdim(L)=(0|3)$}: $L={\bf L_{0|3}^1}=\left<0|e_1,e_2,e_3\right>$.\\
    \item \underline{$\sdim(L)=(1|2)$};  $L=\left<e_1|e_2,e_3\right>$.

    \begin{multicols}{2}
        
        \begin{enumerate}
            \item ${\bf L_{1|2}^1} =\left<e_1|e_2,e_3\right>$ \textup{(}abelian\textup{)};

            \item ${\bf L_{1|2}^2}=\left<e_1|e_2,e_3; [e_2,e_3]=e_1 \right>$;

        \columnbreak        
            \item ${\bf L_{1|2}^3}=\left<e_1|e_2,e_3;[e_1,e_2]=e_3\right>$;

            \item ${\bf L_{1|2}^4}=\left<e_1|e_2,e_3;[e_3,e_3]=e_1\right>$.
         
        \end{enumerate}
    \end{multicols}    
    \item \underline{$\sdim(L)=(2|1)$}: $L=\left<e_1, e_2|e_3\right>$.
        \begin{multicols}{2}
        \begin{enumerate}
            \item ${\bf L_{2|1}^1} =\left<e_1,e_2|e_3\right>$ \textup{(}abelian\textup{)};

        \columnbreak             
            \item ${\bf L_{2|1}^2} =\left<e_1,e_2|e_3;[e_3,e_3]=e_2\right>$.
           
        \end{enumerate}
        \end{multicols}
    \item \underline{$\sdim(L)=(3|0)$}: $L=\left<e_1, e_2,e_3\right>$.
        \begin{multicols}{2}
        \begin{enumerate}
            \item ${\bf L_{3|0}^1} =\left<e_1,e_2,e_3\right>$ \textup{(}abelian\textup{)};

        \columnbreak             
            \item ${\bf L_{3|0}^2} =\left<e_1,e_2,e_3; [e_1,e_2]=e_3\right>$. 
          
        \end{enumerate}
        \end{multicols}
\end{Lemma}~

\begin{proof} This proof is composed of several cases. \\

    \noindent\underline{The case where $\sdim(L)=(0|3)$.} Here, $\dim(L_\ev)=0$ and all the brackets vanish. The only nilpotent Lie superalgebra of super-dimension $\sdim(L)=(0|3)$ is abelian.\\

    \noindent\underline{The case where $\sdim(L)=(1|2)$.} The even part $L_\ev$ is spanned by, let us say,  $e_1$. By adding two odd elements $e_2$ and $e_3$, we obtain the following multiplication table:   
    $$[e_1,e_2]=\lambda_2e_2+\lambda_3e_3,~[e_2,e_2]=\gamma_1e_1,~[e_1,e_3]=\mu_2e_2+\mu_3e_3,~[e_2,e_3]=\gamma_2e_1,~[e_3,e_3]=\gamma_3e_1,$$ with $\lambda_2,\lambda_3,\gamma_1,\gamma_2,\gamma_3,\mu_2,\mu_3\in\K.$ 
    It follows that $\lambda_2=\mu_3=0$, otherwise the Lie superalgebra is not nilpotent. Using the Jacobi identity, we get an algebraic system that leads to the following classes of 
 Lie superalgebras:

    \begin{enumerate}
    \item[(i)]
    ${\bf L^1(\gamma_2,\gamma_3)}$ : $\gamma_1=1,\lambda_3=\mu_2=0,\gamma_2,\gamma_3$ arbitrary.
    \item[(ii)] ${\bf L^2(\gamma_2,\gamma_3)}$ : $\gamma_1=0,\gamma_2=1,\lambda_3=\mu_2=0,\gamma_2\neq 0,\gamma_3$ arbitrary.
    \item[(iii)] ${\bf L^3(\lambda_3,\mu_2)}$ : $\gamma_1=\gamma_2=\gamma_3=0,\lambda_3,\mu_2$ arbitrary.
    \item[(iv)] ${\bf L^4(\gamma_3)}$ : $\gamma_1=\gamma_2=0,\lambda_3=\mu_2=0, \gamma_3$ arbitrary.
    \end{enumerate}
    
    \vspace{0.2cm}

    \begin{itemize}
        \item Study of ${\bf L^4(\gamma_3)}$. First, let us suppose that $\gamma_3=0$. It follows that ${\bf L^4(0)}={\bf L_{1|2}^1}$. let us suppose now  that $\gamma_3\neq 0.$ It follows that ${\bf L^4(\gamma_3)}$ is isomorphic to ${\bf L^4(1)}={\bf L_{1|2}^4}$. 
        \item Study of ${\bf L^3(\lambda_3,\mu_2)}$. If $\lambda_3$ and $\mu_2$ are both non zero, then ${\bf L^3(\lambda_3,\mu_2)}$ is not  nilpotent. Thus, we have three cases to consider. 
             \begin{enumerate}
                 \item Case $\lambda_3=\mu_2=0$: in this case ${\bf L^3(0,0)}={\bf L_{1|2}^1}$.
                 \item Case $\lambda_3\neq 0,\mu_2=0$: in this case ${\bf L^3(\lambda_3,0)}$ is isomorphic to ${\bf L^3(1,0)}={\bf L_{1|2}^3}$.
                 \item  Case $\lambda_3= 0,\mu_2\neq0$: in this case ${\bf L^3(0,\mu_3)}$ is isomorphic to ${\bf L^3(0,1)}$, which is isomorphic to ${\bf L_{1|2}^3}$ as well.
           \end{enumerate}
        \item Study of ${\bf L^2(\gamma_2,\gamma_3)}$. The superalgebra ${\bf L^2(\gamma_2,\gamma_3)}$ is isomorphic to ${\bf L^2(1,0)}={\bf L_{1|2}^2}.$
        \item Study of ${\bf L^1(\gamma_2,\gamma_3)}$. There are four cases to consider.
        \begin{enumerate}
            \item Case $\gamma_2=\gamma_3=0.$ The Lie superalgebra ${\bf L^1(0,0)}$ is then abelian and therefore isomorphic to ${\bf L_{1|2}^1}$.
            \item Case $\gamma_2\neq0,\gamma_3=0.$ The Lie superalgebra ${\bf L^1(\gamma_2,0)}$ is isomorphic to ${\bf L^1(1,0)}$ which in turn is isomorphic to ${\bf L_{1|2}^2}$. Explicitly, the isomorphism $\varphi:{\bf L_{1|2}^2}\longrightarrow {\bf L^1(1,0)}$ is given by   $$\varphi(e_1)=e_1,~\varphi(e_2)=e_3,~\varphi(e_3)=e_2-\frac{1}{2}e_3. $$
            \item Case $\gamma_2=0,\gamma_3\neq 0.$ The Lie superalgebra ${\bf L^1(0,\gamma_3)}$ is isomorphic to ${\bf L^1(0,1)}$, which in turn is isomorphic to ${\bf L_{1|2}^2}$ as well. Explicitly, the isomorphism $\varphi:{\bf L_{1|2}^2}\longrightarrow {\bf L^1(0,1)}$ is given by
            $$\varphi(e_1)=-2ie_1,~\varphi(e_2)=-ie_2+e_3,~\varphi(e_3)=e_2-ie_3, $$ where $i=\sqrt{-1}$.
            \item Case $\gamma_2\neq 0,\gamma_3\neq 0.$ Let us denote $\Gamma:=\gamma_2^2-\gamma_3$. There are two sub-cases, depending on whether $\Gamma=0$ or $\Gamma\neq 0$.
            \begin{enumerate}
                \item Sub-case $\Gamma=0$. In that case, the Lie superalgebra ${\bf L^1(\gamma_2,\gamma_2^2)}$ is isomorphic to ${\bf L_{1|2}^4}$ and the isomorphism $\varphi:{\bf L_{1|2}^4}\longrightarrow {\bf L^1(\gamma_2,\gamma_2^2)}$ is given by
            $$\varphi(e_1)=\gamma_2^2e_1,~\varphi(e_2)=-\gamma_2e_2+e_3,~\varphi(e_3)=e_3. $$

                \item Sub-case $\Gamma\neq0,~\Gamma=\sqrt{\Gamma}$. In that case, the Lie superalgebra ${\bf L^1(\gamma_2,\gamma_3)}$ is isomorphic to ${\bf L_{1|2}^2}$ and the isomorphism $\varphi:{\bf L_{1|2}^2}\longrightarrow {\bf L^1(\gamma_2,\gamma_3)}$ is given by
            $$\varphi(e_1)=-2e_1,~\varphi(e_2)=-\left(1+\gamma_2\right)e_2+e_3,~\varphi(e_3)=\left(1-\gamma_2\right)e_2+e_3. $$
            
                \item Sub-case $\Gamma\neq \sqrt{\Gamma}$. In that case, the superalgebra ${\bf L^1(\gamma_2,\gamma_3)}$ is isomorphic to ${\bf L_{1|2}^2}$ as well and the isomorphism $\varphi:{\bf L_{1|2}^2}\longrightarrow {\bf L^1(\gamma_2,\gamma_3)}$ is given by
            $$\varphi(e_1)=\left(\Gamma-\sqrt{\Gamma}\right)e_1,~\varphi(e_2)=-\left(\gamma_2-\sqrt{\Gamma}\right)e_2+e_3,~\varphi(e_3)=\left(\sqrt{\Gamma}-\gamma_2\right)e_2+e_3. $$
            
            \end{enumerate}
        \end{enumerate}
    \end{itemize}
    As a result, we obtained four nilpotent Lie superalgebras of super-dimension $\sdim(L)=(1|2)$, which are denoted by ${\bf L_{1|2}^1},{\bf L_{1|2}^2},{\bf L_{1|2}^3},{\bf L_{1|2}^4}$ as above.\\

\noindent\underline{Case $\sdim(L)=(2|1)$.} Denote by $e_1$ and $e_2$ the basis element of the  $2$-dimensional even part $L_\ev$. Observe that $[L_\ev, L_\ev] =0$ as $L_\ev$ is nilpotent. As we add one odd element $e_3$, we get the following multiplication table:
    $$[e_1,e_3]=\lambda_1e_3,~[e_2,e_3]=\lambda_2e_3,~[e_3,e_3]=\mu_1e_1+\mu_2e_2,$$ with $\lambda_1,\lambda_2,\mu_1,\mu_2\in\K.$ Observe that  $\lambda_1=\lambda_2=0$, otherwise the Lie superalgebra is not nilpotent. Thus, there are two cases to consider.

    \begin{enumerate}
        \item Case $\mu_1=\mu_2=0$. In this case, the Lie superalgebra is abelian and is denoted by ${\bf L_{2|1}^1}$.
        \item Case $(\mu_1,\mu_2)\neq (0,0)$. In that case, the Lie  superalgebra is isomorphic to ${\bf L_{2|1}^2}$, whose bracket is given by $[e_3,e_3]=e_2$. Explicitly, the isomorphism $\varphi:{\bf L}\longrightarrow {\bf L_{2|1}^2}$ is given by
            $$
            \varphi(e_1)=\left \{\begin{array}{ll}-\mu_2e_1+\frac{1-\mu_2^2}{\mu_1}e_2, & \text{if } \mu_1 \not =0 \\
            e_1, & \text{if } \mu_1 =0 
            \end{array}\right. ,~\varphi(e_2)=\left\{
            \begin{array}{ll}
            \mu_1e_1+\mu_2e_2, & \text{if } \mu_1 \not =0\\
            \frac{1}{\mu_2}e_2,& \text{if } \mu_1 =0
        \end{array} \right.,  ~\varphi(e_3)=e_3.
        $$
    \end{enumerate}
    
\noindent\underline{Case $\sdim(L)=(3|0)$.} This case is covered by \cite{SU}.\\
\end{proof}

With our non-isomorphic Lie superalgebras now found, we can investigate $p|2p$-structures. Here is the strategy. Let $L=L_\ev\oplus L_\od$ be a restricted Lie superalgebra. The even part $L_\ev$ is a restricted Lie algebra and therefore it is equipped with a $p$-map, denoted $(\cdot)^{[p]}$. Using Jacobson's Theorem \ref{SJac}, the $p$-map $(\cdot)^{[p]}$ defines a $p|2p$ structure on $L$ if and only if  
\begin{equation}\label{eqSJac}
    \bigl[e_i^{[p]},f_j\bigl]=\ad^p_{e_i}(f_j),~\forall e_i \text{ basis elements of } L_\ev,~\forall f_j \text{ basis elements of } L_\od.
\end{equation}
Therefore, it is sufficient to check the identity (\ref{eqSJac}) on basis elements of $L_\od$ knowing that the $p$-structures on $L_\ev$ have been classified in \cite{SU}. 

\ssbegin{Theorem}\label{classif3}
Let $L$ be a $p$-nilpotent Lie superalgebra of total dimension $3$ over an
algebraically closed field $\K$ of characteristic $p\geq3$. The  equivalence classes of the $p|2p$-maps on $L$ are given by

\begin{itemize}
    \item \underline{$\sdim(L)=(0|3)$}: $L=\left<0|e_1,e_2,e_3\right>$: none.\\
    \item \underline{$\sdim(L)=(1|2)$}: $L=\left<e_1|e_2,e_3\right>$.

    \begin{multicols}{2}
        
        \begin{enumerate}
            \item ${\bf L_{1|2}^1} =\left<e_1|e_2,e_3\right>$ \textup{(}abelian\textup{)}:
                    \begin{enumerate}
                        \item $e_1^{[p]}=0;$
                    \end{enumerate}
            \item ${\bf L_{1|2}^2}=\left<e_1|e_2,e_3; [e_2,e_3]=e_1 \right>$:
                    \begin{enumerate}
                        \item $e_1^{[p]}=0;$
                    \end{enumerate}
        \columnbreak        
            \item ${\bf L_{1|2}^3}=\left<e_1|e_2,e_3;[e_1,e_2]=e_3\right>$:
                    \begin{enumerate}
                        \item $e_1^{[p]}=0.$
                    \end{enumerate}
            \item ${\bf L_{1|2}^4}=\left<e_1|e_2,e_3;[e_3,e_3]=e_1\right>$:
                    \begin{enumerate}
                        \item $e_1^{[p]}=0;$
                    \end{enumerate}            
        \end{enumerate}
    \end{multicols}    
    \item \underline{$\sdim(L)=(2|1)$}: $L=\left<e_1, e_2|e_3\right>$.
        \begin{multicols}{2}
        \begin{enumerate}
            \item ${\bf L_{2|1}^1} =\left<e_1,e_2|e_3\right>$ \textup{(}abelian\textup{)}:
                    \begin{enumerate}
                        \item $e_1^{[p]}=e_2^{[p]}=0;$
                        \item $e_1^{[p]}=e_2,~e_2^{[p]}=0.$
                    \end{enumerate}
        \columnbreak             
            \item ${\bf L_{2|1}^2} =\left<e_1,e_2|e_3;[e_3,e_3]=e_2\right>$:
                    \begin{enumerate}
                        \item $e_1^{[p]}=e_2^{[p]}=0;$
                        \item $e_1^{[p]}=e_2,~e_2^{[p]}=0.$
                    \end{enumerate}            
        \end{enumerate}
        \end{multicols}
    \item \underline{$\sdim(L)=(3|0)$}: $L=\left<e_1, e_2,e_3\right>$, see \cite{SU}.
        \begin{multicols}{2}
        \begin{enumerate}
            \item ${\bf L_{3|0}^1} =\left<e_1,e_2,e_3\right>$ \textup{(}abelian\textup{)}:
                    \begin{enumerate}
                        \item $e_1^{[p]}=e_2^{[p]}=e_3^{[p]}=0;$
                        \item $e_1^{[p]}=e_2,~e_2^{[p]}=e_3^{[p]}=0;$
                        \item $e_1^{[p]}=e_2,~e_2^{[p]}=e_3,~e_3^{[p]}=0.$
                    \end{enumerate}
        \columnbreak             
            \item ${\bf L_{3|0}^2} =\left<e_1,e_2,e_3; [e_1,e_2]=e_3\right>$ 
                    \begin{enumerate}
                        \item $e_1^{[p]}=e_2^{[p]}=e_3^{[p]}=0;$
                        \item $e_1^{[p]}=e_3,~e_2^{[p]}=e_3^{[p]}=0$.\\
                    \end{enumerate}            
        \end{enumerate}
        \end{multicols}
\end{itemize}

\end{Theorem}

\begin{proof} Cases to study are listed here. ~\\

\noindent\underline{Case $\sdim(L)=(0|3)$.} Since $\dim(L_\ev)=0$, there is no $p$-structure.\\

\noindent\underline{Case $\sdim(L)=(1|2)$.} There is only one 
$p$-nilpotent $p$-structure on $L_\ev$, given by $e_1^{[p]}=0$, which always satisfies Equation \eqref{eqSJac}.\\

\noindent\underline{Case $\sdim(L)=(2|1)$.} As proven in \cite{SU}, there are two non-isomorphic $p$-structures on the two-dimensional Lie algebra $L_\ev$, given by $e_1^{[p]_1}=e_2^{[p]_1}=0$ and $e_1^{[p]_2}=e_2,~e_2^{[p]_2}=0$.
\begin{itemize}
    \item Case of ${\bf L^1_{2|1}}$. As the Lie superalgebra is abelian, Equation (\ref{eqSJac}) is satisfied and both $(\cdot)^{[p]_1}$ and $(\cdot)^{[p]_2}$ define a $p|2p$ structure on ${\bf L^1_{2|1}}$.
     \item Case of ${\bf L^2_{2|1}}$. We have $\ad_{e_1}(e_3)=\ad_{e_2}(e_3)=0$, therefore Equation (\ref{eqSJac}) is satisfied and both $(\cdot)^{[p]_1}$ and $(\cdot)^{[p]_2}$ define a $p|2p$ structure on ${\bf L^2_{2|1}}$.
\end{itemize}

\noindent\underline{Case $\sdim(L)=(3|0)$.} For this case, we refer to \cite{SU}.\end{proof}

\section{$p$-nilpotent restricted  Lie superalgebras of dimension $4$}\label{sectiondim4}


This section aims to classify $p$-nilpotent Lie superalgebras with total dimension 4 over an algebraically closed field of  characteristic $p>2$. First, let's explain our strategy.

\subsection{The strategy} Here we present our method for classifying $p$-nilpotent Lie superalgebras of dimension $4$ adapted from \cite{DU,MS}. We will show that any $n$-dimensional nilpotent Lie superalgebra can be obtained as a central extension of a $(n-1)$-dimensional nilpotent Lie superalgebras by means of a suitable $2$-cocycle. The details of the construction are given below. 

Let $L$ be a restricted Lie superalgebra. We denote by $\Aut(L)$ (resp. $\Aut_p(L)$) the automorphism group (resp. the restricted automorphism group) of $L$. Let $M$ be a restricted $L$-module. The restricted $2$-cocycles space $Z^2_*(L;M)$ is endowed with an $\Aut_p(L)$-structure in the following way. Let $(\varphi,\omega)\in Z^2_*(L;M)$ and $A\in \Aut_p(L)$. The action of $\Aut_p(L)$ on $Z^2_*(L;M)$ is given by $A\cdot(\varphi,\omega):=(A\varphi,A\omega)$, with

\begin{equation}\label{action}
    \begin{cases}
        (A\varphi)(x,y)&=\varphi\left(A(x),A(y)\right),~\forall x,y\in L;\\
        (A\omega)(x)&=\omega\left(A(x)\right),~\forall x\in L_\ev.      
    \end{cases}
\end{equation}

One can show that the restricted $2$-coboundaries space $B^2_*(L;M)$ is an $\Aut_p(L)$-submodule of $Z^2_*(L;M)$ under this action, therefore it is possible to endow the quotient space $H^2_*(L;M)$ with an $\Aut_p(L)$-module structure. We say that two cocycles are \textit{equivalent} if they belong to the same $\Aut_p(L)$-orbit. We show that if two elements of $H^2_*(L;M)$ belong to the same $\Aut_p(L)$-orbit, the central extensions obtained are isomorphic. The action \eqref{action} clearly induces an action of $\Aut(L)$ on $Z^2_{\ce}(L;M)$.  The first step is to determine the $\Aut(L)$-orbits of every Lie superalgebra given in Theorem \ref{classif3} by picking up a representative and to build the associated central extension to compute new Lie brackets. However, as explained in \cite{MS}, the converse is not true and thus we have to  remove isomorphic Lie superalgebras from the list we obtained. At this point, it remains to investigate the $p|2p$-structures on every superalgebra we obtained, which can be done in the same fashion as in the proof of Theorem \ref{classif3}, using the classification of \cite{SU} and Jacobson's Theorem \ref{SJac}.\\

\noindent\textbf{Remark.} In the following, we will only consider Lie superalgebras $L=L_\ev\oplus L_\od$ of total dimension $3$ such that $\dim(L_\od)\geq 1$ to build Lie superalgebras of total dimension $4$. The purely even case has been considered in \cite{SU} using a different approach. We will not miss any case by doing so.


\subsection{Useful results} The following section present some important theoretical results used in our classification scheme for restricted Lie superalgebras of dimension $4$.

\sssbegin{Proposition}
Let $\bigl(L,[\cdot,\cdot]\bigl)$ be a nilpotent Lie superalgebra and let $\varphi$ and 
$\psi$ be two equivalent scalar $2$-cocycles. Let $\bigl(L_{\varphi},[\cdot,\cdot]_{\varphi}\bigl)$ and $\bigl(L_{\psi},[\cdot,\cdot]_{\psi}\bigl)$ be two central extensions of $L$ by an element $X$ such that $|X|=|\varphi|=|\psi|$, with brackets given by
$$[\cdot,\cdot]_{\varphi}:=[\cdot,\cdot]+\varphi(\cdot,\cdot)X $$ and $$[\cdot,\cdot]_{\psi}:=[\cdot,\cdot]+\psi(\cdot,\cdot)X  $$ respectively. Then, the Lie superalgebras $L_{\varphi} $ and $L_{\psi} $ are isomorphic.
\end{Proposition}

\begin{proof} Since $\varphi$ and $\psi$ are equivalent cocycles, there exists $A\in\Aut(L)$ such that $$\varphi(x,y)=\psi(Ax,Ay)~\forall x,y\in L.$$ Consider the map $f:L_{\varphi}\rightarrow L_{\psi}$ given by $f(x)=Ax~\forall x\in L,~f(X)=X.$ Clearly, $f$ is 
 a bijective linear map. Let $x,y\in L$.
\begin{align*}
    f\bigl([x,y]_{\varphi}\bigl)&=A\bigl([x,y]\bigl)+\varphi(x,y)f(X)\\
                                &=[Ax,Ay]+\psi\bigl(Ax,Ay\bigl)X\\
                                &=[f(x),f(y)]_{\psi}.
\end{align*}
Therefore, the map $f$ is an isomorphism and $L_{\varphi}$ and $L_{\psi}$ are isomorphic as Lie superalgebras.
\end{proof}

\sssbegin{Lemma}\label{restrictedh}
    Let $\left( L ,[\cdot,\cdot]_{ L },(\cdot)^{[p]_{ L }}\right)$ be a $p$-nilpotent restricted Lie superalgebra and let $ I\subset \fz(L)$ be a purely odd central ideal of $L$. Then, the quotient space $ H= L / I$ is a restricted Lie superalgebra.
\end{Lemma}

\begin{proof}
  We define a $p$-map $(\cdot)^{[p]_{ H}}$ on $ H$ by $\bar{x}^{[p]_{ H}}:=x^{[p]_{ L }}~\forall \bar{x}\in H_\ev$. This map is well defined. Indeed, let $\bar{x},\bar{y}\in H_\ev$ such that $\bar{x}=\bar{y}$. Then, $x-y\in I$ and therefore $x=y$, because $ I$ is purely odd.
  
  Let $\bar{x}\in H_\ev$ and $\bar{y}\in H$. We have
        $$\bigl[\bar{x}^{[p]_{ H}},\bar{y}\bigl]_{ H}=\bigl[\bar{x}^{[p]_{ L }},y+ I\bigl]_{ L }+ I=\ad^p_x(y)+ I=\ad^p_{\bar{x}}(\bar{y}).$$ Since $ I\subset \fz(L)$, we have $s_j(x,y)=s_j(\bar{x},\bar{y})~\forall 1\leq j\leq p-1.$ For $\bar{x},\bar{y}\in H_\ev$, we have
        $$(\bar{x}+\bar{y})^{[p]_{ H}}=(x+y)^{[p]_{ L }}=x^{[p]_{ L }}+y^{[p]_{ L }}+\displaystyle\sum^{p-1}_{j=1}s_j(x,y)=\bar{x}^{[p]_{ H}}+\bar{y}^{[p]_{ H}}+\displaystyle\sum^{p-1}_{j=1}s_j(\bar{x},\bar{y}).$$ We also have $(\lambda\bar{x})^{[p]_{ H}}=\lambda^p\bar{x}^{[p]_{ H}},~\forall \lambda\in\K$. Therefore the map $(\cdot)^{[p]_{ H}}$ defines a $p|2p$-structure on $ H$.
\end{proof}

\sssbegin{Proposition}\label{propext1dim}
Let $ L $ be a $p$-nilpotent restricted Lie superalgebra of dimension $n$. Then, $ L $ is isomorphic to a central extension by a restricted $2$-cocycle of a $p$-nilpotent restricted Lie superalgebra of dimension $n-1$.
\end{Proposition}

\begin{proof}
    Since $ L $ is $p$-nilpotent, we have $\fz(L)\neq 0$. There are two cases to consider, depending on whether the center admits an even element or not.\\
    
    \underline{Case $\fz(L)_{\ev}\neq 0$}. Since $ L $ is $p$-nilpotent, there is a central $p$-ideal $ I\subset \fz(L)_{\ev}$ such that $ I^{[p]}=0$. In that case, the result follows from \cite[Theorem 2.7]{DU}\footnote{Note that the published version of the paper omits the proof, but it can be found on the arXiv version (\texttt{arXiv:1412.8377}, Lemma 2.2.2).}. In that case, the restricted $2$-cocycle $(\varphi,\omega)$ used to build the extension belongs to $Z^2_*(L;I)^+_\ev$ (see Lemma \ref{shapecocycles}).  \\   

    \underline{Case $\fz(L)_{\ev}=\{0\}$}. In that case, the center is purely odd. Let $Y\in \fz(L)$ and let $ I$ be the one-dimensional ideal generated by the odd element $Y$. We have $Y^{[2p]}=\left(\frac{1}{2}[Y,Y]\right)^{[p]}=0$, so $ I$ is a $p|2p$-ideal. Moreover, with Lemma \ref{restrictedh} we have that the quotient space $ H= L / I$ is a restricted Lie superalgebra. We have the short exact sequence of restricted Lie superalgebras
$$0\longrightarrow I\longrightarrow L \overset{\pi}{\longrightarrow} H\longrightarrow 0,   $$
with $\pi:  L \rightarrow H$ the projection map, which is a restricted Lie superalgebras morphism.
Then, we choose a (linear) section $\rho: H\rightarrow L $ of $\pi$, which is an even linear map satisfying $\pi\circ\rho=\id_{ H}$. We define an even bilinear map $\phi: H\times H\rightarrow I$ by $\phi(x,y)=[\rho(x),\rho(y)]-\rho([x,y])$. One can show that $\text{Im}(\phi)\subset\Ker(\pi)= I$. As a consequence, there exists an odd bilinear form $\Delta: H\times H\rightarrow\K$ such that $\phi(x,y)=\Delta(x,y)Y,~\forall x,y\in H$.\\

\textbf{Claim:} $\phi\in Z^2_{\ce}( H; I)_\ev$ and  $\Delta\in Z^2_{\ce}( H;\K)_\od$.\\
To prove this claim, let $x,y,z\in H$. We have 
    \begin{align*}
        d^2_{\ce}\phi(x,y,z)&=(-1)^{|x||z|}\phi(x,[y,z])+(-1)^{|x||y|}\phi(y,[z,x])+(-1)^{|y||z|}\phi(z,[x,y])\\
        &=(-1)^{|x||z|}\Bigl(\bigl[\rho(x),\rho([y,z])\bigl]-\rho\bigl([x,[y,z]]\bigl)\Bigl)
        +(-1)^{|x||y|}\Bigl(\bigl[\rho(y),\rho([z,x])\bigl]-\rho\bigl([y,[z,x]]\bigl)\Bigl)\\
        &+(-1)^{|y||z|}\Bigl(\bigl[\rho(z),\rho([x,y])\bigl]-\rho\bigl([z,[x,y]]\bigl)\Bigl)\\
        &=(-1)^{|x||z|}\Bigl[\rho(x),\bigl[\rho(y),\rho(z)\bigl]-\phi(y,z)\Bigl]+(-1)^{|x||y|}\Bigl[\rho(y),\bigl[\rho(z),\rho(x)\bigl]-\phi(z,x)\Bigl]\\
        &+(-1)^{|y||z|}\Bigl[\rho(z),\bigl[\rho(x),\rho(y)\bigl]-\phi(x,y)\Bigl]\\
        &-\underset{=~0}{\rho\Bigl(\underbrace{(-1)^{|x||z|}\bigl([x,[y,z]]\bigl)+(-1)^{|x||y|}\bigl([y,[z,x]]\bigl)+(-1)^{|y||z|}\rho\bigl([z,[x,y]]\bigl)}\Bigl)}\\
        &=(-1)^{|x||z|}\Bigl[\rho(x),\bigl[\rho(y),\rho(z)\bigl]\Bigl]+(-1)^{|x||y|}\Bigl[\rho(y),\bigl[\rho(z),\rho(x)\bigl]\Bigl]\\
        &+(-1)^{|y||z|}\Bigl[\rho(z),\bigl[\rho(x),\rho(y)\bigl]\Bigl]\\
        &=0.
    \end{align*}
    Therefore, $\phi\in Z^2_{\ce}( H; I)_\ev$ and an immediate consequence is that $\Delta\in Z^2_{\ce}( H;\K)_\od$.
    
    We are now ready to prove that $ L $ is isomorphic to the one-dimensional central extension of $ H$ by the odd cocycle $\Delta$, that we will denote $ H_{\Delta}$. Let $x\in L $ and set $z=x-\rho\circ\pi(x)$. Then, $\pi(z)=0$ and $z\in I$. Let $y=\pi(x)\in H$. We therefore have the decomposition $x=\rho(y)+z$. We consider the bijective linear map $f: L \rightarrow H_{\Delta}$ given by $f(x)=y+z$. It remains to show that $f$ is a restricted Lie morphism. Let $x_1,x_2\in L $. There exists $y_1,y_2\in H$ and $z_1,z_2\in I$ such that $x_1=\rho(y_1)+z_1$ and $x_2=\rho(y_2)+z_2$. 
We have 
$$f\bigl([x_1,x_2]\bigl)=f\Bigl(\bigl[\rho(y_1)+z_1,\rho(y_2)+z_2\bigl]\Bigl)=f\Bigl(\phi(y_1,y_2)+\rho\bigl([y_1,y_2]\bigl)\Bigl)=\phi(y_1,y_2)+[y_1,y_2].$$ On the other hand, we have 
$$\bigl[f(x_1),f(x_2)\bigl]=\Bigl[f\bigl(\rho(y_1)+z_1,\rho(y_2)+z_2\bigl)\Bigl]=\bigl[y_1+z_1,y_2+z_2\bigl]=\phi(y_1,y_2)+[y_1,y_2].$$ The map $f$ is a Lie morphism. For the restricted part, let $x\in L _\ev$. The decomposition is then given by $x=\rho(y)$, since $x$ is even. We have
$$f\left(x^{[p]}\right)=f\left(\rho(y)^{[p]}\right)=f\left(\rho\bigl(y^{[p]}\bigl)\right)=y^{[p]}=f(x)^{[p]}.$$ Therefore, $f$ is a restricted Lie isomorphism and the restricted superalgebra $ L $ is obtained as an one-dimensional central extension of $ H$ by means of a restricted cocycle of the form $(\Delta,0)$, with $\Delta\in Z^2_{\ce}( H; \K)_\od$.

\end{proof}

\subsection{Equivalent $2$-cocycles}

 The first step in our strategy is to determine the
orbits of $2$-cocycles for every Lie superalgebra listed in Theorem \ref{classif3}. Two restricted $2$-cocycles $(\varphi,\omega)$ and $(\psi,\theta)$ are said to be \textit{equivalent} if there exists $A\in\Aut_p(L)$ such that $(\varphi,\omega)=A\cdot(\psi,\theta)$, that is to say, they are in the same orbit with respect to the action (\ref{action}).\\

 Let $L=L_\ev\oplus L_\od$ be a restricted nilpotent Lie superalgebra of superdimension $\sdim(L)=(n|m)$. We denote its basis elements by $e_1,\cdots, e_{n+m}$ with
 $$L_\ev=\Span(e_1,\cdots,e_n)~\text{and}~L_\od=\Span(e_{n+1},\cdots,e_{n+m}).$$
 A basis for (ordinary) $2$-cocycles is then given by elements of the form
\begin{equation}\label{deltacocycles}\Delta_{i,j}:L\times L\longrightarrow \K,~1\leq i\leq n+m,~i\leq j\leq n+m,\end{equation}
 where $\Delta_{i,j}(e_k,e_l)=\delta_{i,k}\delta_{j,l}$ and $\Delta_{i,j}=-(-1)^{|e_i||e_j|}\Delta_{j,i}$.

 Our method allows us to consider only ordinary $2$-cocycles, the ``restricted part" of the structure being obtained by direct computation using Jacobson's Theorem \ref{SJac}.

 For every Lie superalgebra listed in the following Lemmas, we first compute the list of $2$-cocycles using the Mathematica package \texttt{SuperLie} (see \cite{G}). We then remove from the list those that  are in the same orbit.

   Suppose now that $L$ is of total dimension $3$. An element $\tau\in \mathfrak{S}_3$ induces an automorphism $A_{\tau}$ of $L$ given by
    \begin{equation}
        A_{\tau}(e_i)=e_{\tau^{-1}(i)}.
    \end{equation}
    Therefore, we deduce an action of $\Aut(L)$ on $Z_{\ce}^2(L;\K)$ given by
    \begin{equation}\label{action3}   A_{\tau}\cdot\Delta_{i,j}=\Delta_{\tau(i),\tau(j)},~\tau\in\mathfrak{S}_3,~\Delta_{i,j}\in Z^2_{\ce}(L;M).  
    \end{equation}   
     For every superalgebra, we first give the conditions on $A$ so that it is a Lie morphism, which needs to be completed with the condition $\det(A)\neq 0$ (that we will not always write but that we must always check). We will use the notation $\Delta\cong\Delta'$ if $\Delta$ and $\Delta'$ are equivalent $2$-cocycles.

    \sssbegin{Lemma}
     Suppose that $L$ is a nilpotent Lie superalgebra with $\sdim(L)=(0|3)$ over an algebraically closed field of characteristic $p\geq3$. The equivalence classes of (ordinary) homogeneous $2$-cocycles on $L$ are given by
     $0,~\Delta_{1,1},~~\Delta_{1,2},~\Delta_{1,1}+\Delta_{2,3}.$ 
    \end{Lemma}

\begin{proof}
     The general form of an automorphism $A\in\Aut(L)$ is given by (where  $\lambda_i,\mu_i,\gamma_i,\in\K$ for $i=1,2,3$):
\begin{equation}\label{autoauto}
Ae_1=\lambda_1e_1+\lambda_2e_2+\lambda_3e_3,~Ae_2=\mu_1 e_1+\mu_2 e_2+\mu_3 e_3,~Ae_3=\gamma_1 e_1+\gamma_2 e_2+\gamma_3 e_3.
\end{equation}
We have $\dim\left(H^2_{\ce}(L;\K)\right)=6$ and the basis is given by 
$$\left\lbrace\Delta_{1,1},~\Delta_{2,2},~\Delta_{3,3},~\Delta_{1,2},~\Delta_{1,3},~\Delta_{2,3}\right\rbrace.$$

Using the action \eqref{action3}, we have $\Delta_{1,1}\cong\Delta_{2,2}\cong\Delta_{3,3}$ and $\Delta_{1,2}\cong\Delta_{1,3}\cong\Delta_{2,3}$. Moreover, for all non-zero $a\in\K$, we have $a\Delta_{i,j}\cong\Delta_{i,j},~\forall 1\leq i,j\leq 3$. Then, we consider cocycles of size $2$, \textit{i.e.} of the form $a\Delta_{i,j}+b\Delta_{k,l},~a,b\in \K^{\times}$. The list can be reduced using the action \eqref{action3}. We show that every cocycle of this form is equivalent to $\Delta_{1,1},~\Delta_{1,2}$ or $\Delta_{1,1}+\Delta_{2,3}$. Next, we consider cocycles of size $3,4,5,6$. Each time, we reduce the list using \eqref{action3} and we show that every cocycle is equivalent to $\Delta_{1,1},~\Delta_{1,2}$ or $\Delta_{1,1}+\Delta_{2,3}$. The method is technical and consists in choosing an automorphism of the form \eqref{autoauto} and computing the appropriate values of the parameters using the first part of Equation \eqref{action} on the basis. Since the computations are lengthy and include numerous case-by-case analyses, we omit them.
\end{proof}
    
  \sssbegin{Lemma}
     Suppose that $L$ is a nilpotent Lie superalgebra with $\sdim(L)=(1|2)$ over an algebraically closed field of characteristic $p\geq3$. The equivalence classes of (ordinary) $2$-cocycles on $L$ are given by
    
    \begin{itemize}
        \item[] $L={\bf L^1_{1|2}}$: $0,~\Delta_{1,2},~\Delta_{2,3},~\Delta_{2,2}+\Delta_{2,3}+\Delta_{3,3};$
        \item[] $L={\bf L^2_{1|2}}$: $0,~\Delta_{2,2},~\Delta_{2,2}+\Delta_{3,3};$
        \item[] $L={\bf L^3_{1|2}}$: $0,~\Delta_{1,3},~\Delta_{2,2};$
        \item[] $L={\bf L^4_{1|2}}$: $0,~\Delta_{2,2},~\Delta_{2,3},~\Delta_{2,2}+\Delta_{2,3}.$
    \end{itemize}
 \end{Lemma}

 \begin{proof}
      The general form of an automorphism $A\in\Aut(L)$ is given by
$$Ae_1=\lambda e_1,~Ae_2=\mu_2 e_2+\mu_3 e_3,~Ae_3=\gamma_2 e_2+\gamma_3 e_3,~\lambda,\mu_2,\mu_3,\gamma_2,\gamma_3\in\K.$$ The condition on 
$\det(A)$ is given by $\lambda(\mu_2\gamma_3-\mu_3\gamma_2)\neq 0$.

    \begin{itemize}

        \item $L={\bf L^1_{1|2}}$. No conditions on $A$, except $\det(A)\neq 0$. We have $\dim\left(H^2_{\ce}(L;\K)\right)=5$ and there are two odds cocycles, namely $\Delta_{1,2} $ and $\Delta_{1,3}$; and three even cocycles, namely $\Delta_{2,2}$, $\Delta_{3,3}$ and $\Delta_{2,3}$ in the basis.
        For all non-zero $a,b\in\K$, we have $a\Delta_{1,2}\cong\Delta_{1,2}\cong\Delta_{1,3}\cong b\Delta_{1,3}.$ Moreover, we have $\Delta_{1,2}\cong a\Delta_{1,2}+b\Delta_{1,3},$ the automorphism being given by
        $$ Ae_1=e_1,~Ae_2=ae_2,~ Ae_3=be_2.      $$ We therefore have one non-trivial orbit of odd cocycles represented by $\Delta_{1,2}$.
        For the even cocycles, we first have $a\Delta_{2,3}\cong\Delta_{2,3},~\forall a\in\K$. Then, for all non-zero $a,b\in\K$, we have $\Delta_{2,3}\cong a\Delta_{2,2}+b\Delta_{3,3}$ with the automorphism given by
        $$Ae_1=e_1,~Ae_2=e_2+i\sqrt{\frac{a}{b}}e_3,~Ae_3=\frac{1}{2a}e_2-\frac{i}{2\sqrt{ab}}e_3,~\text{with}~i=\sqrt{-1}.$$ Moreover, for all non-zero $a,b\in\K$, we have $\Delta_{2,3}\cong a\Delta_{2,3}+b\Delta_{2,2}$ with the automorphism given by
        $$Ae_1=e_1,~Ae_2=e_2+\frac{b}{2}e_3,~Ae_3=ae_3.$$
        Finally, let $c,d,e$ be non-zero scalars. It remains to investigate the equivalence classes of the cocycle $c\Delta_{2,3}+d\Delta_{2,2}+e\Delta_{3,3}.$ There are two cases.
            \begin{itemize}
                \item[$(i)$] \underline{Case $c^2\neq de.$} In that case, we have $c\Delta_{2,3}+d\Delta_{2,2}+e\Delta_{3,3}\cong \Delta_{2,3}$ with the automorphism given by
                $$Ae_1=e_1,~Ae_2=\frac{c+\sqrt{c^2-de}}{e}e_2+\frac{2c-e\mu_2}{2}e_3,~Ae_3=e_2+\frac{e}{2}e_3,$$ where $\mu_2=\frac{c+\sqrt{c^2-de}}{e}$. 
                \item[$(ii)$] \underline{Case $c^2=de.$} In that case, we have $c\Delta_{2,3}+d\Delta_{2,2}+e\Delta_{3,3}\cong \Delta_{2,3}+\Delta_{2,2}+\Delta_{3,3}$, the automorphism being given by
                $$Ae_1=e_1,~Ae_2=-\frac{1}{\sqrt{e}}e_3,~Ae_3=e_2-\left(\frac{d}{c}+\frac{1}{\sqrt{e}}\right)e_3.$$                
            \end{itemize}
        Therefore, we have two non-trivial orbits of even cocycles, namely $\Delta_{2,3}$ and $\Delta_{2,3}+\Delta_{2,2}+\Delta_{3,3}$.

        \item $L={\bf L^2_{1|2}}$. The conditions on $A$ are $\lambda=\mu_3\gamma_2+\mu_2\gamma_3,$ and $\mu_2\mu_3=\gamma_2\gamma_3=0$. Since $\lambda$ must be non-zero, we can choose $\mu_2=\gamma_3=0$. We have $\dim\left(H^2_{\ce}(L;\K)\right)=2$ spanned by $\Delta_{2,2}$ and $\Delta_{3,3}$, which are equivalent using the transposition $\tau=(2,3)$. For all non-zero $a,b\in\K$, we have $a\Delta_{2,2}\cong\Delta_{3,3}$ and $b\Delta_{3,3}\cong\Delta_{2,2}$.    Finally, $a\Delta_{2,2}+b\Delta_{3,3}\cong \Delta_{2,2}+\Delta_{3,3}$, the automorphism being given by
        $$Ae_1=\sqrt{ab}e_1,~Ae_2=\sqrt{a}e_3,~Ae_3=\sqrt{b}e_2.  $$        
        Then, one can show that $\Delta_{2,2}+\Delta_{3,3}\ncong\Delta_{2,2}$, that $\Delta_{2,2}+\Delta_{3,3}\ncong 0$ and that $\Delta_{2,2}\ncong 0$. Therefore, we have three orbits: $0,\Delta_{2,2},\Delta_{2,2}+\Delta_{3,3}$.

        \item $L={\bf L^3_{1|2}}$. We have $\dim\left(Z^2_{\ce}(L;\K)\right)=2$ spanned by $\Delta_{2,2}$ and $\Delta_{1,3}$, which are not equivalent since $\Delta_{2,2}$ is even and $\Delta_{1,3}$ is odd.

        \item $L={\bf L^4_{1|2}}$. The conditions on $A$ are $\lambda=\gamma_3^2$ and $\mu_3=0$. We have $\dim\left(Z^2(L;\K)\right)=2$ spanned by $\Delta_{2,2}$ and $\Delta_{2,3}$. First, for all non-zero $a,b\in\K$, we have $a\Delta_{2,2}\cong \Delta_{2,2}$ and $b\Delta_{2,3}\cong \Delta_{2,3}$. However, $\Delta_{2,2}\ncong \Delta_{2,3}$. One can also easily prove that 
        $\Delta_{2,2}\ncong \Delta_{2,2}+\Delta_{2,3} $ and $\Delta_{2,3}\ncong \Delta_{2,2}+\Delta_{2,3}$. It remains to show that $a\Delta_{2,2}+b\Delta_{2,3}\cong \Delta_{2,2}+\Delta_{2,3}$ for all non-zero $a,b\in\K$, which can be achieved with the isomorphism given by
        $$Ae_1=\frac{b^2}{a}e_1,~ Ae_2=\sqrt{a}e_2,~Ae_3=\frac{b}{\sqrt{a}}e_3. \qed $$ 
    \end{itemize}\noqed  \end{proof}

  \sssbegin{Lemma}
     Suppose that $L$ is a nilpotent Lie superalgebra with $\sdim(L)=(2|1)$ over an algebraically closed field of characteristic $p\geq3$. The equivalence classes of \textup{(}ordinary\textup{)} homogeneous $2$-cocycles on $L$ are given by
    
    \begin{itemize}
         \item[] $L={\bf L^1_{2|1}}$: $0,~\Delta_{1,2},~\Delta_{1,3},~\Delta_{3,3},~\Delta_{1,2}+\Delta_{3,3}$;
         \item[] $L={\bf L^2_{2|1}}$: $0,~\Delta_{1,3}$.
    \end{itemize}
 \end{Lemma}

\begin{proof}   
  The general form of an automorphism $A\in\Aut(L)$ is given by
$$Ae_1=\lambda_1e_1+\lambda_2e_2,~Ae_2=\mu_1 e_1+\mu_2 e_2,~Ae_3=\gamma e_3,\text{ where }\lambda_1,\lambda_2,\mu_1,\mu_2,\gamma\in\K.$$ The condition on 
$\det(A)$ is given by $\gamma(\mu_2\lambda_1-\mu_1\lambda_2)\neq 0$.
\begin{itemize}
        \item $L={\bf L^1_{2|1}}$. No conditions on $A$, except $\det(A)\neq 0$. We have $\dim\left(H^2_{\ce}(L;\K)\right)=4$ and there are two even cocycles, namely $\Delta_{1,2},~\Delta_{3,3}$ and two odd cocycles, namely $\Delta_{1,3}$ and $\Delta_{3,3}$ in the basis. For all non-zero $a,b\in\K$, we have $a\Delta_{1,2}\cong\Delta_{1,2}$, $b\Delta_{3,3}\cong\Delta_{3,3}$ but $a\Delta_{1,2}\ncong b\Delta_{3,3}$. Moreover, we have $\Delta_{1,2}+\Delta_{3,3}\cong a\Delta_{1,2}+b\Delta_{3,3}$ with the automorphism given by
        $$Ae_1=ae_1,~Ae_2=e_2,~Ae_3=\sqrt{b}e_3.$$ We therefore have three classes of non-zero even cocycles, namely $\Delta_{1,2},~\Delta_{3,3}$ and $\Delta_{1,2}+\Delta_{3,3}$.

        Furthermore, for all non-zero $a,b\in\K$, we have $ a\Delta_{1,3}\cong b\Delta_{2,3}$ with the automorphism given by 
        $$Ae_1=\frac{b}{a}e_2,~Ae_2=\frac{b}{a}e_1,~Ae_3=e_3.$$ Moreover, $ a\Delta_{1,3}\cong \Delta_{1,3}$ and $ b\Delta_{2,3}\cong \Delta_{2,3}$. Finally, we have $\Delta_{1,3}\cong a\Delta_{1,3}+b\Delta_{2,3}$ with the automorphism given by
        $$Ae_1=\frac{1}{b}e_2,~Ae_2=e_1-\frac{a}{b}e_2,~Ae_3=e_3.$$ There is just one class of non-zero odd cocycles given by $\Delta_{1,3}$.
        
        \item $L={\bf L^2_{2|1}}$. The conditions of $A$ are $\mu_1=0$ and $\gamma^2=\mu_2$. We have $\dim\left(H^2_{\ce}(L;\K)\right)=1$ spanned by $\Delta_{1,3}$, which is not equivalent to zero. Moreover, $a\Delta_{1,3}\cong\Delta_{1,3}$ for all non-zero $a\in\K$.\qed
\end{itemize}
\noqed 
\end{proof}
The following Theorem summarizes the results.

 \sssbegin{Theorem}\label{cocycles}
     Suppose that $L$ is a nilpotent Lie superalgebra of total dimension $3$ with $\dim(L_\od)\geq 1$ over an algebraically closed field of characteristic $p\geq3$. The equivalence classes of \textup{(}ordinary\textup{)}  homogeneous $2$-cocycles on $L$ are given by
    
    \begin{itemize}
        \item[] $L={\bf L^1_{0|3}}$:  $0,~\Delta_{1,1},~~\Delta_{1,2},~\Delta_{1,1}+\Delta_{2,3};$ 
        \item[] $L={\bf L^1_{1|2}}$: $0,~\Delta_{1,2},~\Delta_{2,3},~\Delta_{2,2}+\Delta_{2,3}+\Delta_{3,3};$
        \item[] $L={\bf L^2_{1|2}}$: $0,~\Delta_{2,2},~\Delta_{2,2}+\Delta_{3,3};$
        \item[] $L={\bf L^3_{1|2}}$: $0,~\Delta_{1,3},~\Delta_{2,2};$
        \item[] $L={\bf L^4_{1|2}}$: $0,~\Delta_{2,2},~\Delta_{2,3},~\Delta_{2,2}+\Delta_{2,3}.$
        \item[] $L={\bf L^1_{2|1}}$: $0,~\Delta_{1,3},~\Delta_{1,2},~\Delta_{3,3},~\Delta_{1,2}+\Delta_{3,3}$;
        \item[] $L={\bf L^2_{2|1}}$: $0,~\Delta_{1,3}$.
    \end{itemize}
 \end{Theorem}

\subsection{Building the extensions}\label{extensions}

Following are tables presenting Lie superalgebras obtained through the central extension process using non-equivalent cocycles of Theorem  \ref{cocycles}. To each Lie superalgebra of total dimension $4$ obtained this way, we give the scalar cocycle used to construct it, the parity of the added element $X$ as well as the new Lie bracket. We use the following notation:
\begin{equation}
  [x,y]_{\text{new}}=[x,y]_{\text{old}}+
  \Delta(x,y) X.
\end{equation}
As explained in Prop. \eqref{propext1dim}, an extension of a Lie superalgebra $L$ by an odd ideal of dimension $1$ can be seen as a extension of $L$ by the scalar field using an odd cocycle $\Delta$. Therefore, the degree of the added element $X$ is by definition $|X|=|\Delta|$.

{\renewcommand{\arraystretch}{1.5} 

\begin{center}
    \begin{tabular}{|c|c|c|c|c|}
    \hline
          \textbf{Name} &\textbf{sdim} & \textbf{Cocycle}& \textbf{Added element} &\textbf{Bracket} \\\hline
          ${\bf L^a_{2|2}}$ & $(2|2)$ & $0$ & $X$ even & $[\cdot,\cdot]=0$ \\\hline
         ${\bf L^a_{1|3}}$ &$(1|3)$ & $0$ & $X$ odd & $[\cdot,\cdot]=0$  \\\hline
         ${\bf L^b_{1|3}}$& $(1|3)$  & $\Delta_{1,2}$   & $X$ odd & $[e_1,e_3]=X$\\\hline
         ${\bf L^b_{2|2}}$& $(2|2)$  & $\Delta_{2,3}$   & $X$ even & $[e_2,e_3]=X$\\\hline
         ${\bf L^c_{2|2}}$& $(2|2)$  & $\Delta_{2,2}+\Delta_{2,3}+\Delta_{3,3}$   & $X$ even & $[e_2,e_2]=[e_2,e_3]=[e_3,e_3]=X$\\\hline
    \end{tabular}
   \\~\\ Lie superalgebras obtained by extensions of ${\bf L^1_{1|2}}$.
\end{center}

\vspace{0.5cm}

\begin{center}
    \begin{tabular}{|c|c|c|c|c|}
    \hline
          \textbf{Name} &\textbf{sdim} & \textbf{Cocycle}& \textbf{Added element} &\textbf{Bracket} \\\hline
          ${\bf L^d_{2|2}}$ & $(2|2)$ & $0$ & $X$ even & $[e_2,e_3]=e_1$ \\\hline
         ${\bf L^c_{1|3}}$ &$(1|3)$ & $0$ & $X$ odd & $[e_2,e_3]=e_1$  \\\hline
         ${\bf L^e_{2|2}}$& $(2|2)$  & $\Delta_{2,2}$   & $X$ even & $[e_2,e_3]=e_1,~[e_2,e_2]=X$\\\hline
         ${\bf L^f_{2|2}}$& $(2|2)$  & $\Delta_{2,2}+\Delta_{3,3}$   & $X$ even & $[e_2,e_3]=e_1,~[e_2,e_2]=[e_3,e_3]=X$\\\hline
    \end{tabular}
   \\~\\ Lie superalgebras obtained by extensions of ${\bf L^2_{1|2}}$.
\end{center}
    
\vspace{0.5cm}

\begin{center}
    \begin{tabular}{|c|c|c|c|c|}
    \hline
          \textbf{Name} &\textbf{sdim} & \textbf{Cocycle}& \textbf{Added element} &\textbf{Bracket} \\\hline
          ${\bf L^g_{2|2}}$ & $(2|2)$ & $0$ & $X$ even & $[e_1,e_2]=e_3$ \\\hline
         ${\bf L^d_{1|3}}$ &$(1|3)$ & $0$ & $X$ odd & $[e_1,e_2]=e_3$   \\\hline
         ${\bf L^e_{1|3}}$ &$(1|3)$ & $\Delta_{1,3}$ & $X$ odd & $[e_1,e_2]=e_3,~[e_1,e_3]=X$   \\\hline
         ${\bf L^h_{2|2}}$& $(2|2)$  & $\Delta_{2,2}$   & $X$ even & $[e_1,e_2]=e_3,~[e_2,e_2]=X$\\\hline
    \end{tabular}
   \\~\\ Lie superalgebras obtained by extensions of ${\bf L^3_{1|2}}$.
\end{center}
    
\vspace{0.5cm}

\begin{center}
    \begin{tabular}{|c|c|c|c|c|}
    \hline
          \textbf{Name} &\textbf{sdim} & \textbf{Cocycle}& \textbf{Added element} &\textbf{Bracket} \\\hline
          ${\bf L^i_{2|2}}$ & $(2|2)$ & $0$ & $X$ even & $[e_3,e_3]=e_1$ \\\hline
         ${\bf L^f_{1|3}}$ &$(1|3)$ & $0$ & $X$ odd & $[e_3,e_3]=e_1$   \\\hline
         ${\bf L^j_{2|2}}$& $(2|2)$  & $\Delta_{2,2}$   & $X$ even & $[e_3,e_3]=e_1,~[e_2,e_2]=X$\\\hline
         ${\bf L^k_{2|2}}$& $(2|2)$  & $\Delta_{2,3}$   & $X$ even & $[e_3,e_3]=e_1,~[e_2,e_3]=X$\\\hline
         ${\bf L^l_{2|2}}$& $(2|2)$  & $\Delta_{2,3}+\Delta_{2,3}$   & $X$ even & $[e_3,e_3]=e_1,~[e_2,e_3]=[e_2,e_3]=X$\\\hline
    \end{tabular}
   \\~\\ Lie superalgebras obtained by extensions of ${\bf L^4_{1|2}}$.
\end{center}
    
\vspace{0.5cm}

\begin{center}
    \begin{tabular}{|c|c|c|c|c|}
    \hline
          \textbf{Name} &\textbf{sdim} & \textbf{Cocycle}& \textbf{Added element} &\textbf{Bracket} \\\hline
          ${\bf L^a_{3|1}}$ & $(3|1)$ & $0$ & $X$ even & $[\cdot,\cdot]=0$ \\\hline
         ${\bf L^m_{2|2}}$ &$(2|2)$ & $0$ & $X$ odd & $[\cdot,\cdot]=0$   \\\hline
         ${\bf L^n_{2|2}}$& $(2|2)$  & $\Delta_{1,3}$   & $X$ odd  & $[e_1,e_3]=X$\\\hline
         ${\bf L^b_{3|1}}$& $(3|1)$  & $\Delta_{1,2}$   & $X$ even & $[e_1,e_2]=X$\\\hline
         ${\bf L^c_{3|1}}$& $(3|1)$  & $\Delta_{3,3}$   & $X$ even & $[e_3,e_3]=X$\\\hline
         ${\bf L^d_{3|1}}$& $(3|1)$  & $\Delta_{1,2}+\Delta_{3,3}$   & $X$ even & $[e_1,e_2]=[e_3,e_3]=X$\\\hline
    \end{tabular}
   \\~\\ Lie superalgebras obtained by extensions of ${\bf L^1_{2|1}}$.
\end{center}

\vspace{0.5cm}

\begin{center}
    \begin{tabular}{|c|c|c|c|c|}
    \hline
          \textbf{Name} &\textbf{sdim} & \textbf{Cocycle}& \textbf{Added element} &\textbf{Bracket} \\\hline
          ${\bf L^e_{3|1}}$ & $(3|1)$ & $0$ & $X$ even & $[e_3,e_3]=e_2$ \\\hline
         ${\bf L^o_{2|2}}$ &$(2|2)$ & $0$ & $X$ odd & $[e_3,e_3]=e_2$   \\\hline
         ${\bf L^p_{2|2}}$& $(2|2)$  & $\Delta_{1,3}$   & $X$ odd  & $[e_3,e_3]=e_2,~ [e_1,e_3]=X$\\\hline         
    \end{tabular}
   \\~\\ Lie superalgebras obtained by extensions of ${\bf L^2_{2|1}}$.
\end{center}

\vspace{0.5cm}

\begin{center}
    \begin{tabular}{|c|c|c|c|c|}
    \hline
          \textbf{Name} &\textbf{sdim} & \textbf{Cocycle}& \textbf{Added element} &\textbf{Bracket} \\\hline
          ${\bf L^a_{0|4}}$ & $(0|4)$ & $0$ & $X$ odd & $[\cdot,\cdot]=0$ \\\hline
         ${\bf L^g_{1|3}}$ &$(1|3)$ & $0$ & $X$ even & $[\cdot,\cdot]=0$  \\\hline
         ${\bf L^h_{1|3}}$& $(1|3)$  & $\Delta_{1,1}$   & $X$ even & $[e_1,e_1]=X$\\\hline
         ${\bf L^i_{1|3}}$& $(1|3)$  & $\Delta_{1,2}$   & $X$ even & $[e_1,e_2]=X$\\\hline
         ${\bf L^j_{1|3}}$& $(1|3)$  & $\Delta_{1,1}+\Delta_{2,3}$  & $X$ even & $[e_1,e_2]=[e_2,e_3]=X$\\\hline
    \end{tabular}
   \\~\\ Lie superalgebras obtained by extensions of ${\bf L^1_{0|3}}$.
\end{center}

\vspace{0.5cm}
{\renewcommand{\arraystretch}{1}

\subsection{Detecting isomorphisms}

Some of the Lie superalgebras in the tables of Section \ref{extensions}  are redundant. In the next step of our procedure, we will remove isomorphic Lie superalgebras from the tables. By simply rearranging and relabeling the basis elements, we are able to obtain the following isomorphisms:\\

\noindent\underline{$\sdim(L)=(1|3)$}: ${\bf L_{1|3}^a}\cong {\bf L_{1|3}^g},$ ${\bf L_{1|3}^b}\cong {\bf L_{1|3}^d},$ ${\bf L_{1|3}^c}\cong {\bf L_{1|3}^i},$ ${\bf L_{1|3}^f}\cong {\bf L_{1|3}^h}.$\\

\noindent\underline{$\sdim(L)=(2|2)$}: ${\bf L_{2|2}^a}\cong {\bf L_{2|2}^m},$ ${\bf L_{2|2}^b}\cong {\bf L_{2|2}^d},$ ${\bf L_{2|2}^e}\cong {\bf L_{2|2}^k},$ ${\bf L_{2|2}^g}\cong {\bf L_{2|2}^n},$ ${\bf L_{2|2}^h}\cong {\bf L_{2|2}^p},$ ${\bf L_{2|2}^i}\cong {\bf L_{2|2}^o}.$\\

\noindent\underline{$\sdim(L)=(3|1)$}: ${\bf L_{3|1}^c}\cong {\bf L_{3|1}^e}.$ \\

\subsubsection{Tables of invariants}
In order to classify the remaining Lie superalgebras, we compute some  classical invariants, namely the derived superalgebra, the superdimension of the center, and the superdimension of the cohomology spaces, which have been calculated using the Mathematica package \texttt{SuperLie} (see \cite{G}) for $p=3,5,7,11$. The dimension may vary depending on the characteristic, especially for $p=3$. We indicate it in the tables.\\
{\renewcommand{\arraystretch}{1.5}
\footnotesize{
\begin{center}
    \begin{tabular}{|c|c|c|c|c|c|c|}
    \hline
          $L$ & $[L,L]$ & $\sdim(\fz(L))$ & $\sdim\left(H^1_{\ce}(L;\K)\right)$& $\sdim\left(H^2_{\ce}(L;\K)\right)$ & $\sdim\left(H^3_{\ce}(L;\K)\right)$&$\sdim\left(H^4_{\ce}(L;\K)\right)$ \\\hline
          ${\bf L^a_{1|3}}$   &  $0$ & $1|3$  & $1|3$  &  $6|3$ & $7|9$  & $15|10$ \\\hline
          ${\bf L^b_{1|3}}$   & $\left\langle X  \right\rangle$  &  $0|2$ & $ 1|2 $  &$ 3|2 $    & $ 3|4~(3|5\text{ if}~p=3) $  & $ 5|4~(7|5\text{ if}~p=3) $ \\\hline
          ${\bf L^c_{1|3}}$   & $\left\langle e_1 \right\rangle$  &  $1|1$ &  $ 0|3 $  & $ 5|0 $   & $ 0|7 $   &$ 9|0 $ \\\hline
          ${\bf L^e_{1|3}}$   & $\left\langle e_3,X  \right\rangle$  & $0|1$  &  $ 1|1 $  & $ 2|1 $   & $ 2|2~(2|4\text{ if}~p=3) $   &$ 3|2~(5|4\text{ if}~p=3) $\\\hline
          ${\bf L^f_{1|3}}$   & $\left\langle e_1 \right\rangle$  &  $1|2$ &  $ 0|3 $  & $ 5|0 $   & $ 0|7 $   &$ 9|0 $ \\\hline
          ${\bf L^j_{1|3}}$   & $\left\langle X \right\rangle$  & $1|0$  &  $ 0|3 $  & $ 5|0 $   & $ 0|7 $   &$ 9|0 $\\\hline              
    \end{tabular}
   \\~\\ Invariants for Lie superalgebras of 
 $\sdim=(1|3)$.
\end{center}}
\normalsize{}
\vspace{0.2cm}

\footnotesize{
\begin{center}
    \begin{tabular}{|c|c|c|c|c|c|c|}
    \hline
          $L$ & $[L,L]$ & $\sdim(\fz(L))$ & $\sdim\left(H^1_{\ce}(L;\K)\right)$& $\sdim\left(H^2_{\ce}(L;\K)\right)$ & $\sdim\left(H^3_{\ce}(L;\K)\right)$&$\sdim\left(H^4_{\ce}(L;\K)\right)$ \\\hline
          ${\bf L^a_{2|2}}$   &  $0$ & $2|2$   & $2|2$  &  $4|4$ & $6|6$  & $8|8$ \\\hline
          ${\bf L^b_{2|2}}$   & $\left\langle X \right\rangle$  &  $2|0$ & $ 1|2 $  &$ 2|2 $    & $ 2|2 $   & $ 2|2 $ \\\hline
          ${\bf L^c_{2|2}}$   & $\left\langle X \right\rangle$  &  $2|1$ &  $ 1|2 $  & $ 2|2 $   & $ 2|2 $   & $ 2|2 $ \\\hline
          ${\bf L^e_{2|2}}$   & $\left\langle e_1,X \right\rangle$  & $2|0$  &  $ 0|2 $  & $ 1|1 $   & $ 1|1 $   &$ 1|1 $\\\hline
          ${\bf L^f_{2|2}}$   & $\left\langle e_1,X\right\rangle $ & $2|0$ &  $ 0|2 $  & $ 1|0 $   & $ 0 $   &$ 0 $\\\hline
          ${\bf L^g_{2|2}}$   & $\left\langle e_3 \right\rangle$  &  $1|1$ & $ 2|1 $  &$ 2|2 $    & $ 2|2~(2|3\text{ if}~p=3) $   & $ 2|2~(3|4\text{ if}~p=3) $ \\\hline
          ${\bf L^h_{2|2}}$   & $\left\langle X,e_3 \right\rangle$  &  $1|1$ &  $ 1|1 $  & $ 1|1 $   & $ 1|1~(1|2 \text{ if}~p=3) $   & $ 1|1~(1|2 \text{ if}~p=3) $ \black\\\hline
          ${\bf L^i_{2|2}}$   & $\left\langle e_1 \right\rangle$  & $2|1$  &  $ 1|2 $  & $ 2|2 $   & $ 2|2 $   & $ 2|2 $\\\hline
          ${\bf L^j_{2|2}}$   & $\left\langle e_1,X\right\rangle$ & $2|0$  &  $ 0|2 $  & $ 1|0 $   & $ 0 $   &$ 0 $\\\hline
          ${\bf L^l_{2|2}}$   & $\left\langle e_1,X \right\rangle$  &  $2|0$ & $ 0|2 $  & $ 1|0 $    & $ 0 $   & $ 0 $ \\\hline            
    \end{tabular}
   \\~\\ Invariants for Lie superalgebras of $\sdim=(2|2)$.
\end{center}}
\normalsize{}
\vspace{0.2cm}

\footnotesize{
\begin{center}
    \begin{tabular}{|c|c|c|c|c|c|c|}
    \hline
          $L$ & $[L,L]$ & $\sdim(\fz(L))$ & $\sdim\left(H^1_{\ce}(L;\K)\right)$& $\sdim\left(H^2_{\ce}(L;\K)\right)$ & $\sdim\left(H^3_{\ce}(L;\K)\right)$&$\sdim\left(H^4_{\ce}(L;\K)\right)$ \\\hline
          ${\bf L^a_{3|1}}$   &  $0$ & $3|1$   & $3|1$  &  $4|3$ & $4|4$  & $4|4$\\\hline
          ${\bf L^b_{3|1}}$   & $\left\langle X \right\rangle$  &  $1|1$ & $ 2|1 $  &$ 3|2 $    & $ 3|3 $   & $ 3|3 $ \\\hline
          ${\bf L^c_{3|1}}$   & $\left\langle X \right\rangle$  &  $3|0$ &  $ 2|1 $  & $ 1|2 $   & $ 0|1 $   &$ 0 $ \\\hline
          ${\bf L^d_{3|1}}$   & $\left\langle X \right\rangle$  & $1|0$  &  $ 2|1 $  & $ 1|2 $   & $ 0|1 $   &$ 0 $\\\hline          
    \end{tabular}
   \\~\\ Invariants for Lie superalgebras of $\sdim=(3|1)$.
\end{center}}
\normalsize{}
\vspace{0.5cm}

\sssbegin{Lemma}~
\begin{itemize}
    \item \underline{$\sdim(L)=(1|3)$}. The Lie superalgebras ${\bf L_{1|3}^a}$, ${\bf L_{1|3}^b}$, ${\bf L_{1|3}^c}$, ${\bf L_{1|3}^e}$, ${\bf L_{1|3}^f}$, ${\bf L_{1|3}^j}$ are pairwise non-isomorphic. 
    \item \underline{$\sdim(L)=(2|2)$}. The Lie superalgebras  ${\bf L_{3|1}^a}$, ${\bf L_{3|1}^b}$, ${\bf L_{3|1}^e}$, ${\bf L_{3|1}^f}$, ${\bf L_{3|1}^g}$, ${\bf L_{3|1}^h}$, ${\bf L_{3|1}^i}$ are pairwise non-isomorphic. Moreover, we have ${\bf L_{2|2}^f}\cong {\bf L_{2|2}^j},$ ${\bf L_{2|2}^f}\cong {\bf L_{2|2}^l}$ and ${\bf L_{2|2}^c}\cong {\bf L_{2|2}^i}$.
    \item \underline{$\sdim(L)=(3|1)$}. The Lie superalgebras ${\bf L_{3|1}^a}$, ${\bf L_{3|1}^b}$, ${\bf L_{3|1}^c}$, ${\bf L_{3|1}^d}$ are pairwise non-isomorphic.
\end{itemize}

\end{Lemma}
\begin{proof}
The ``non-isomorphic" results of the Lemma follow from the tables of invariants. As for the rest, we exhibit explicit isomorphisms.\\
     An isomorphism $\phi:{\bf L_{2|2}^j}\rightarrow {\bf L_{2|2}^f}$ is given by
     $$\phi(e_1)=2X-2e_1,~\phi(X)=2X+2e_1,~\phi(e_2)=e_2+e_3,~\phi(e_3)=e_2-e_3.$$
     An isomorphism $\varphi:{\bf L_{2|2}^f}\rightarrow {\bf L_{2|2}^l}$ is given by
     $$\varphi(e_1)=X-e_1,~\varphi(X)=e_1,~\varphi(e_2)=e_3,~\varphi(e_3)=e_2-e_3.$$
     An isomorphism $\psi:{\bf L_{2|2}^c}\rightarrow {\bf L_{2|2}^i}$ is given by
     $$\psi(e_1)=X+e_1,~\psi(X)=e_1,~\psi(e_2)=e_2+e_3,~\psi(e_3)=e_3.$$
\end{proof}

Hereafter, we summarize the classification of $4$-dimensional Lie superalgebras. We have relabeled the basis elements and have renamed the Lie superalgebras. For clarity, we indicate the ``working name" of the Lie superalgebras.

\sssbegin{Theorem}\label{classif4} The classification of $4$-dimensional nilpotent Lie superalgebras over an algebraically closed field of characteristic different from $2$ is given by:\\
\noindent\underline{$\sdim(L)=(0|4)$}: $L=\left\langle 0|x_1,x_2,x_3,x_4 \right\rangle$

${\bf L_{0|4}^1}:~[\cdot,\cdot]=0.$

\noindent\underline{$\sdim(L)=(1|3)$}: $L=\left\langle x_1|x_2,x_3,x_4 \right\rangle$

${\bf L_{1|3}^1}~(={\bf L_{1|3}^a}):$ abelian;

${\bf L_{1|3}^2}~(={\bf L_{1|3}^b}):~[x_1,x_3]=x_4;$

${\bf L_{1|3}^3}~(={\bf L_{1|3}^c}):~[x_2,x_3]=x_1;$

${\bf L_{1|3}^4}~(={\bf L_{1|3}^e}):~[x_1,x_2]=x_3,~[x_1,x_3]=x_4;$

${\bf L_{1|3}^5}~(={\bf L_{1|3}^f}):~[x_3,x_3]=x_1;$

${\bf L_{1|3}^6}~(={\bf L_{1|3}^j}):~[x_2,x_2]=x_1,~[x_3,x_4]=x_1.$

\noindent\underline{$\sdim(L)=(2|2)$}: $L=\left\langle x_1,x_2|x_3,x_4 \right\rangle$

${\bf L_{2|2}^1}~(={\bf L_{2|2}^a}):$ abelian;

${\bf L_{2|2}^2}~(={\bf L_{2|2}^b}):~[x_3,x_4]=x_2;$

${\bf L_{2|2}^3}~(={\bf L_{2|2}^e}):~[x_3,x_3]=x_2,~[x_3,x_4]=x_1;$

${\bf L_{2|2}^4}~(={\bf L_{2|2}^f}):~[x_3,x_3]=[x_4,x_4]=x_2,~[x_3,x_4]=x_1;$

${\bf L_{2|2}^5}~(={\bf L_{2|2}^g}):~[x_1,x_3]=x_4;$

${\bf L_{2|2}^6}~(={\bf L_{2|2}^h}):~[x_1,x_3]=x_4,~[x_3,x_3]=x_2.$

${\bf L_{2|2}^7}~(={\bf L_{2|2}^i}):~[x_4,x_4]=x_1.$

\noindent\underline{$\sdim(L)=(3|1)$}: $L=\left\langle x_1,x_2,x_3|x_4 \right\rangle$

${\bf L_{3|1}^1}~(={\bf L_{3|1}^a}):$ abelian;

${\bf L_{3|1}^2}~(={\bf L_{3|1}^b}):~[x_1,x_2]=x_3;$

${\bf L_{3|1}^3}~(={\bf L_{3|1}^c}):~[x_2,x_2]=x_3;$

${\bf L_{3|1}^4}~(={\bf L_{3|1}^d}):~[x_1,x_2]=[x_3,x_4]=x_3.$

\noindent\underline{$\sdim(L)=(4|0)$}: $L=\left\langle x_1,x_2,x_3,x_4|0 \right\rangle$ \textup{(}see \cite{SU}\textup{)}

${\bf L_{4|0}^1}:$ abelian;

${\bf L_{4|0}^2}:~[x_1,x_2]=x_3;$

${\bf L_{4|0}^3}:~~[x_1,x_2]=x_3,~[x_1,x_3]=x_4.$
\end{Theorem}

\noindent\textbf{Remark.} These results are in agreement with those in \cite{AH,AH2} over the field of complex numbers.

\subsection{$p|2p$-structures}
We use Jacobson's Theorem \ref{SJac} to investigate $p$-nilpotent structures on the Lie superalgebras listed in the Theorem \ref{classif4}.  Let $L$ be a $p$-nilpotent restricted  Lie superalgebra of dimension $4$. The even part $L_\ev$ is also a $p$-nilpotent Lie algebra of dimension $\leq 4$. The classification of those Lie superalgebras has been achieved in \cite{SU}. Therefore, given a $p$-map $(\cdot)^{[p]}$ on $L_\ev$, it is enough to check the condition
\begin{equation}\label{eqSJac'}
    \bigl[e_i^{[p]},f_j\bigl]=\ad^p_{e_i}(f_j),~\forall e_i \text{ basis elements of } L_\ev,~\forall f_j \text{ basis elements of } L_\od,
\end{equation} to decide whether $(\cdot)^{[p]}$ defines a $p$-nilpotent structure on $L$ or not. With this method, we obtain all $p$-nilpotent restricted Lie superalgebras.

\sssbegin{Theorem}\label{pmap4}
    The $p$-nilpotent structures on nilpotent Lie superalgebras of total dimension $4$ are given by:

    \begin{enumerate}
        \item $\sdim(L)=(0|4)$: none.
        \item $\sdim(L)=(1|3)$: $x_1^{[p]}=0$.
        \item $\sdim(L)=(2|2)$:
            \begin{enumerate}
                \item $x_1^{[p]_1}=x_2^{[p]_1}=0;$
                \item $x_1^{[p]_2}=x_2,~x_2^{[p]_2}=0.$  
            \end{enumerate}        
        \item $\sdim(L)=(3|1)$:
             \begin{enumerate}
                \item Case $L_\ev$ abelian:
                    \begin{enumerate}
                        \item $x_1^{[p]_1}=x_2^{[p]_1}=x_3^{[p]_1}=0;$
                        \item $x_1^{[p]_2}=x_2,~x_2^{[p]_2}=x_3^{[p]_2}=0.$
                        \item $x_1^{[p]_3}=x_2,~x_2^{[p]_3}=x_3,~x_3^{[p]_3}=0.$
                    \end{enumerate}
                \item Case $L_\ev\cong {\bf L_{3|0}^2}$ \textup{(}see Theorem \ref{classif3}\textup{)}:  
                    \begin{enumerate}
                        \item $x_1^{[p]_4}=x_2^{[p]_4}=x_3^{[p]_4}=0;$
                        \item $x_1^{[p]_5}=x_3,~x_2^{[p]_5}=x_3^{[p]_5}=0.$
                    \end{enumerate}
            \end{enumerate}  
        \item $\sdim(L)=(4|0)$: see \cite[Theorem 2.1]{SU}.
    \end{enumerate}
\end{Theorem}

\begin{proof} We will examine a few cases: \\
    \underline{The case where  $\sdim(L)=(1|3)$.} We have $L_\ev=\left\langle x_1\right\rangle$ and the only $p$-nilpotent $p$-map on $L_\ev$ is given by $x_1^{[p]}=0$. Then, $\ad_{x_1}^3=0$ for all Lie superalgebra of superdimension $\sdim(L)=(1|3)$ listed in Theorem \ref{classif4}. As a consequence, Eq. \eqref{eqSJac'} is always satisfied for $p>2$.\\
    \underline{The case where  $\sdim(L)=(2|2)$.} We have $L_\ev=\left\langle x_1, x_2\right\rangle$. There are two $p$-nilpotent $p$-maps given by $x_1^{[p]_1}=x_2^{[p]_1}=0$ and $x_1^{[p]_2}=x_2,~x_2^{[p]_2}=0.$ If $x_1$ and $x_2$ are central, $(\cdot)^{[p]_1}$ always satisfies Eq. \eqref{eqSJac'}. Therefore, it remains to investigate ${\bf L_{2|2}^5}$ and ${\bf L_{2|2}^6}$. For both Lie superalgebras, we have $\ad_{x_2}=0$ and $\ad_{x_1}^2=0$, so Eq. \eqref{eqSJac'} is satisfied. For $(\cdot)^{[p]_2}$, notice that the image of $(\cdot)^{[p]_2}$  always lies in the center and that $\ad^2=0$, which implies that Eq. \eqref{eqSJac'} is also satisfied.\\
    \underline{The case where  $\sdim(L)=(3|1)$.} On the same lines as the previous cases.
\end{proof}

\subsection{Comparison with the purely cohomological methods}\label{comparison} The authors in \cite{MS} employed the restricted cohomology to obtain $5$-dimensional $p$-nilpotent Lie algebras as restricted central extensions. To build $p$-maps on the extensions, they used equivalence classes of \textit{restricted} $2$-cocycles. In the present work, we did not use the ``restricted part" of the $2$-cocycles, since the possible $p|2p$ structures are completely determined on the even part of the Lie superalgebra. The purpose of this section is to compare our method on an example with the ``purely cohomological" method described in \cite{MS}.\\

Consider $L={\bf L_{1|2}^2}=\left\langle e_1|e_2,e_3,~[e_2,e_3]=e_1,~e_1^{[p]}=0\right\rangle$. The non-equivalent non-zero (ordinary) $2$-cocycles are $\phi_1=\Delta_{2,2}$ and $\phi_2=\Delta_{2,2}+\Delta_{3,3}$. For $j\in\{1,2\}$, a $p$-homogeneous map $\omega:L_\ev\rightarrow\K$ is $\phi_j$-compatible if
\begin{equation}
\label{star39}
\omega(x+y)=\omega(x)+\omega(y)+\displaystyle\sum_{\underset{x_1=x,~x_2=y}{x_i\in\{x,y\}}}\frac{1}{\sharp(x)}\phi_j\biggl([\cdots[x_1,x_2],x_3],\cdots,x_{p-1}],x_p\biggl),~\forall x,y\in L.
\end{equation}

Suppose that $p>3$. It follows that Eq \eqref{star39} reduces to $\omega(x+y)=\omega(x)+\omega(y)$. Therefore, any $p$-semilinear map $L_\ev\rightarrow\K$ is $\phi_j$-compatible and is of the form $\omega=\alpha e_1^*,~\alpha\in\K$. The ``restricted part" of the cocycle condition is given by
\begin{equation}
\bigl[e_1^{[p]},e_1\bigl]=\phi_j\Bigl(e_1,[e_1,\cdots,[e_1,e_1]\cdots]\Bigl),
\end{equation} must always be satisfied. 

\subsubsection{Case $\phi_1=\Delta_{2,2}$} We compute the equivalences classes of restricted $2$-cocycles. Two restricted cocycles $(\phi_1,\omega)$ and $(\phi_1,\theta)$ are equivalent if there exists a restricted automorphism $A\in\Aut_p(L)$ such that $A\cdot(\phi_1,\theta)=(\phi_1,\omega)$. Therefore, we must have $A\cdot\phi_1=\phi_1$ and $A\cdot\theta=\omega$. Let $\alpha\neq 0$ and consider the restricted automorphism given by
$$A(e_1)=\frac{1}{\alpha^{1/p}}e_1,~A(e_2)=e_2,~A(e_3)=\frac{1}{\alpha^{1/p}}e_3.  $$ Then, we have $A\cdot(\phi_1,\alpha e_1^*)=(\phi_1, e_1^*)$. We therefore have two classes of restricted $2$-cocycles given by $(\phi_1, 0)$ and $(\phi_1, e_1^*)$.
We can build two extensions of ${\bf L_{1|2}^2}$:

\begin{itemize}
    \item ${\bf L^2_{1|2}(\phi_1,0)}=\left\langle e_1,X|e_2,e_3,~ [e_2,e_2]=X,[e_2,e_3]=e_1,~e_1^{[p]}=X^{[p]}=0\right\rangle;$
    \item ${\bf L^2_{1|2}(\phi_1,e_1^*)}=\left\langle e_1,X|e_2,e_3,~ [e_2,e_2]=X,[e_2,e_3]=e_1,~e_1^{[p]}=X,X^{[p]}=0\right\rangle.$
\end{itemize}

We have ${\bf L^2_{1|2}(\phi_1,0)}\cong {\bf L^3_{2|2}}$ with the $p|2p$-map $(3 a)$ (see Theorem \ref{pmap4}) and ${\bf L^2_{1|2}(\phi_1,e_1^*)}\cong {\bf L^3_{2|2}}$ with the $p|2p$-map $(3 b)$ (see Theorem \ref{pmap4}).

\subsubsection{Case $\phi_2=\Delta_{2,2}+\Delta_{3,3}$} In this case, the computations show that there are infinitely many non-equivalent restricted $2$-cocycles, namely $(\phi_2,\alpha e_1^*),~\alpha\in\K$, leading to the following extensions:
    $${\bf L^2_{1|2}(\phi_2,\alpha e_1^*)}=\left\langle e_1,X|e_2,e_3,~[e_2,e_2]=[e_3,e_3]=X,[e_2,e_3]=e_1,~e_1^{[p]_{\alpha}}=\alpha X, X^{[p]}=0\right\rangle.$$

Our goal is to demonstrate that those extensions lead to two non-isomorphic restricted Lie superalgebras. Suppose that $\alpha\neq 0$. We consider the map $A: {\bf L^2_{1|2}(\phi_2,\alpha e_1^*)}\rightarrow {\bf L^2_{1|2}(\phi_2,e_1^*)}$ given by
$$A(e_1)=\sqrt[2p-2]{\alpha}e_1,~A(e_2)=\frac{1}{\alpha}e_2,~A(e_3)=\frac{1}{\alpha}e_4,~A(e_4)=\frac{1}{\alpha}e_3.$$ We claim that this map is a restricted isomorphism of Lie superalgebras. Therefore, $\bf{L^2_{1|2}(\phi_2,\alpha e_1^*)}$ is isomorphic to $\bf{L_{2|2}^5}$ with the $p|2p$-map $(3 b)$ (see Theorem \ref{pmap4}), provided $\alpha\neq 0.$ In the case where $\alpha=0$, we can show that ${\bf L^2_{1|2}(\phi_2,0)}$ is isomorphic to ${\bf L_{2|2}^5}$ with the $p|2p$-map $(3 a)$ (see Theorem \ref{pmap4}).

\section{$p|2p$-structures on the Lie superalgebras $K^{n,m}$, see \cite{GKN}}\label{knm}
The goal of this section is to compute $p|2p$-mappings for various classes of higher dimensional nilpotent Lie superalgebras.

\subsection{The Lie superalgebra $K^{2,m}$, $m$ odd}

The Lie superalgebra $K^{2,m}$ is spanned by the generators $x_0, x_1\;| \; y_1,\ldots y_m$ (Even $|$ Odd), with non-zero brackets given by 
\[
\begin{array}{rlllll}
[x_0, y_i] &= &-[y_i, x_0] &= &y_{i+1}, &  i\leq m-1, \\[2mm] 
[y_i, y_{m+1-i}]& = & [ y_{m+1-i}, y_i] & = &(-1)^{i+1}x_1, & 1 \leq i\leq \frac{m+1}{2}.
\end{array}
\]
It has been shown in \cite{GKN} that a Lie superalgebra of superdimension $n|m$ has a maximal nilindex $n + m-1$ only when $n=2$ and $m$ odd. Moreover, for any odd $m$, there is only one Lie superalgebra with this maximal nilindex and that is $K^{2,m}$. \\

{\bf Claim \cite{BEM}.} $K^{2,m}$ is restricted if and only if $m\leq p$. Explicitly, the $[p|2p]$-map is given by
\[x_0^{[p]}=s_1 x_1, \quad x_1^{[p]}=s_2 x_1, \text{ where $s_1, s_2\in \mathbb K$.}\]

\subsection{The Lie superalgebra $K^{2,m},~m\text{ even}$}

For $m$ even, the Lie superalgebra $K^{2,m}$ is spanned by the generators $x_0, x_1\;| \; y_1,\ldots y_m$ (Even $|$ Odd), with non-zero brackets given by 
\[
\begin{array}{rlllll}
[x_0, y_i] &= &-[y_i, x_0] &= &y_{i+1}, &  1\leq i\leq m-1, \\[2mm] 
[y_i, y_{m-i}]& = & [ y_{m-i}, y_i] & = &(-1)^{(m-2i)/2}x_1, & 1 \leq i\leq \frac{m}{2}.
\end{array}
\]
In that case, it has been shown in \cite{GKN} that the superalgebra $K^{2,m}$ has maximal nilindex $m$.

The superalgebra $K^{2,m},~m\text{ even}$ is restricted if and only if $m\leq p$. On that case, the $[p|2p]$-map is given by
\[x_0^{[p]}=s_1 x_1, \quad x_1^{[p]}=s_2 x_1, \text{ where $s_1, s_2\in \mathbb K$.}\] 

\subsection{The Lie superalgebra $K^{3,m},~m\text{ odd}$}

For $m$ odd, the Lie superalgebra $K^{3,m}$ is spanned by the generators $x_0, x_1,x_2\;| \; y_1,\ldots y_m$ (Even $|$ Odd), with non-zero brackets given by
\[
\begin{array}{rlllll}
[x_0, x_1] &= &-[x_1, x_0] &= &x_{2}, &  \\[2mm] 
[x_0, y_i] &= &-[y_i, x_0] &= &y_{i+1}, &  1\leq i\leq m-1, \\[2mm] 
[y_i, y_{m+1-i}]& = & [ y_{m+1-i}, y_i] & = &(-1)^{(i+1)/2}x_2, & 1 \leq i\leq \frac{m+1}{2}.
\end{array}
\]
We have $K^{3,m}_\ev\cong{\bf L^2_{3|0}}$.
\sssbegin{Proposition}\label{k3mo}
The Lie superalgebra $K^{3,m}$ is restricted if and only if $m\leq p$. In that case, the $p|2p$-map is given by 
\[
x_0^{[p]}=s_0 x_2,~x_1^{[p]}=s_1 x_1,~x_2^{[p]}=s_2 x_2,~s_0,s_1,s_2\in \K.
\]
\end{Proposition}

\subsection{The Lie superalgebra $K^{3,m},~m\text{ even}$}

For $m$ even, the Lie superalgebra $K^{3,m}$ is spanned by the generators $x_0, x_1, x_2\;| \; y_1,\ldots y_m$ (Even $|$ Odd), with non-zero brackets given by 
\[
\begin{array}{rlllll}
[x_0, x_1] &= &-[x_1, x_0] &= &x_{2}, &   \\[2mm] 
[x_0, y_i] &= &-[y_i, x_0] &= &y_{i+1}, &  1\leq i\leq m-1, \\[2mm] 
[y_i, y_{m-i}]& = & [ y_{m-i}, y_i] & = &(-1)^{(m-2i)/2}x_1, & 1 \leq i\leq \frac{m}{2}.\\[2mm] 
[y_i,y_{m+1-i}]& = & [y_{m+1-i},y_i] & = & (-1)^{(m-2i)/2}\left(\frac{1}{2}(m-2i+1)\right)x_2, & 1 \leq i\leq \frac{m}{2}.
\end{array}
\]
In that case, it is not difficult to see that Proposition \ref{k3mo} also holds.

\subsection{The Lie superalgebra $K^{4,m},~m\text{ even}$}

For $m$ even, the Lie superalgebra $K^{4,m}$ is spanned by the generators $x_0, x_1, x_2, x_3\;| \; y_1,\ldots y_m$ (Even $|$ Odd), with non-zero brackets given by 
\[
\begin{array}{rlllll}
[x_0, x_i] &= &-[x_i, x_0] &= &x_{i+1}, &  1\leq i\leq 2. \\[2mm] 
[x_0, y_i] &= &-[y_i, x_0] &= &y_{i+1}, &  1\leq i\leq m-1. \\[2mm] 
[y_i, y_{m-i}]& = & [ y_{m-i}, y_i] & = &(-1)^{(m-2i)/2}x_1, & 1 \leq i\leq \frac{m}{2}.\\[2mm] 
[y_i,y_{m+1-i}]& = & [y_{m+1-i},y_i] & = & (-1)^{(m-2i)/2}\left(\frac{1}{2}(m-2i+1)\right)x_2, & 1 \leq i\leq \frac{m}{2}.\\[2mm] 
[y_i,y_{m+2-i}]& = & [y_{m+2-i},y_i] & = & (-1)^{(m-2i+2)/2}\left(\frac{1}{2}(i-m)(i-1)^2\right)x_3, & 2 \leq i\leq \frac{1}{2}(m+2).
\end{array}
\]
We have $K^{4,m}_\ev\cong{\bf L^3_{4|0}}$.

\sssbegin{Proposition}\label{k4me}
The superalgebra $K^{4,m}$ is restricted if and only if $m<p$. In that case, the $p$-map is given by 
\[
x_0^{[p]}=s_0 x_3,~x_1^{[p]}=s_1 x_3,~x_2^{[p]}=s_2 x_3,~x_3^{[p]}=s_3 x_2,~s_0,s_1,s_2,s_3\in \K.
\]
\end{Proposition}

\subsection{The Lie superalgebra $K^{4,5}$}

The Lie superalgebra $K^{4,5}$ is spanned by the generators $x_0, x_1, x_2, x_3\;| \; y_1,y_2,y_3,y_4,y_5$ (Even $|$ Odd), with non-zero brackets given by 
\[
\begin{array}{rllllll}
[x_0, x_i] &= &-[x_i, x_0] &= &x_{i+1}, & 1\leq i\leq 2& \\ [2mm] 
[x_0, y_i] &= &-[y_i, x_0] &= &x_{i+1}, & 1\leq i\leq 4& \\ [2mm] 
[x_1, y_3] &= &-[y_3, x_1] &= &y_{5}, & &  \\ [2mm] 
[x_2, y_2] &= &-[y_2, x_2] &= &-y_{5}, & &  \\ [2mm]
[x_3, y_1] &= &-[y_1, x_3] &= &y_{5}. & &  \\ [2mm] 
[y_i, y_4] &= &[y_4, y_i] &= &-\frac{3}{2}x_{i+1}, &  1\leq i\leq 2.& \\ [2mm]  
[y_1, y_3] &= &[y_3, y_1] &= &[y_2, y_2]& =~~x_1 &  \\ [2mm] 
[y_2, y_3] &= &[y_3, y_2] &= &\frac{1}{2}x_2& &  \\ [2mm] 
[y_3, y_3] &= & 2x_3 & & & &  
\end{array}
\]
We have $K^{4,5}_\ev\cong{\bf L^3_{4|0}}$.

\sssbegin{Proposition}\label{k45}

For $p=3$, the Lie superalgebra $K^{4,5}$ is restricted with the $3|6$-map given by 
\[
x_0^{[3]}=x_3,~x_1^{[3]}=x_2^{[3]}=x_3^{[3]}=0.
\]

For $p\geq3$, the Lie superalgebra $K^{4,5}$ is restricted with the $p|2p$-map given by 
\[
x_0^{[p]}=x_1^{[p]}=x_2^{[p]}=x_3^{[p]}=0.
\]
\end{Proposition}


\end{document}